\documentclass[a4paper, reqno]{amsart}
\usepackage{amssymb,amsmath,amsfonts,amsthm,mathrsfs}
\usepackage{enumerate}
\usepackage{bbm}
\usepackage{multirow,hhline}
\usepackage{todonotes}

\newcommand{\R}{\mathbb R}
\newcommand{\N}{\mathbb N}
\newcommand{\C}{\mathbb C}
\newcommand{\Z}{\mathbb Z}

\newcommand{\cL}{\mathcal{L}}

\renewcommand{\div}{\mathop{\text{\upshape{div}}}}

\renewcommand{\epsilon}{\varepsilon}
\renewcommand{\bar}[1]{\overline{#1}}

\renewcommand{\tilde}{\widetilde}
\renewcommand{\phi}{\varphi}
\newcommand{\one}{\mathbbm{1}}
\newcommand{\txtd}{\textnormal{d}}

\theoremstyle{definition}

\numberwithin{equation}{section}

\def\Xint#1{\mathchoice
{\XXint\displaystyle\textstyle{#1}}%
{\XXint\textstyle\scriptstyle{#1}}%
{\XXint\scriptstyle\scriptscriptstyle{#1}}%
{\XXint\scriptscriptstyle\scriptscriptstyle{#1}}%
\!\int}
\def\XXint#1#2#3{{\setbox0=\hbox{$#1{#2#3}{\int}$ }
\vcenter{\hbox{$#2#3$ }}\kern-.6\wd0}}

\def\dashint{\Xint-}

\begin{document}

\title[Reduction Methods in Climate Dynamics]
{Reduction Methods in Climate Dynamics\\ - A Brief Review}

\author{Felix Hummel}
\address{Felix Hummel, Department of Mathematics, TU Munich, Germany}

\author{Peter Ashwin}
\address{Peter Ashwin, Department of Mathematics, University of Exeter, Exeter EX4 4QF, UK}

\author{Christian Kuehn}
\address{Christian Kuehn, Department of Mathematics, TU Munich, Germany}

\begin{abstract}
We review a range of reduction methods that have been, or may be useful for connecting models of the Earth's climate system of differing complexity. We particularly focus on methods where rigorous reduction is possible. We aim to highlight the main mathematical ideas of each reduction method and also provide several benchmark examples from climate modelling.
\end{abstract}

\maketitle

\tableofcontents

\section{Introduction}
\label{sec:intro}

In order to understand the Earth's climate system, a larger number of models of varying complexity have been developed: these cover processes at a wide variety of space and timescales and give rise to hierarchy of models \cite{Dijkstra_2013,Ghil_etal_2008} for various aspects of the climate. The highest resolution and most complex Earth System Models (ESMs) such as \cite{CESM2013} that underpin the science Climate Modelling Intercomparison Project (CMIP) forming the scientific basis for the Assessment Reports of the Intergovernmental Panel on Climate Change (IPCC) \cite{IPCC_AR6}. ESMs try to describe the Earth System as accurately as possible and to make precise prediction of the evolution of the Earth's climate and its different components, and from these results we can make conclusions not just about natural variability of the climate system but also of variability that is due to anthropogenic effects. Models that attempt to model the atmosphere and ocean are called general circulation models (GCMs) \cite{Vallis_2017}. Since ESMs and GCMs are of high resolution and include complex physics, they are very difficult to analyze.

Climate and weather models\footnote{We refer to both types of models simply as climate models in the remainder} are typically nonlinear differential equations that are mathematically challenging (compressible Navier-Stokes equations with highly nontrivial thermodynamics and chemistry, and topographically non-trivial boundary conditions and a variety of forcing terms) \cite{Ghil_etal_2008}. Over very long timescales different processes become important for the dynamics become important \cite{ghil2020physics}. Hence, they are usually simulated with a high resolution on very powerful computers but typically only short runs can be achieved at high resolution.  For some questions that involve longer timescales, one can consider a hierarchy of models \cite{Schneider1974} where high resolution ESM/GCMs can be approximated using lower resolution ESMs and simpler Earth system models of intermediate complexity (EMICs). These models usually contain several simplifications and a larger grid size compared to the ESMs/GCMs and they therefore need less computing time. However, EMICs are still mostly too complex for analytical study without further simplification.

To have any chance of analytically understanding components of the Earth's climate, one has to go further down the model hierarchy and study conceptual models such as box models or energy balance models that contain not more than a few variables, where hidden variables may be assumed to be in equilibrium or may be represented as noise terms \cite{Hasselmann_1976,Ghil_etal_2008,Majda_Timofeyev_VandenEijnden_2001,Majda_Timofeyev_VandenEijnden_2003}. It is remarkable that simple models can often reproduce many aspects of the Earth's climate system \cite{Dijkstra_2013,Kaper_Engler_2013} and while the similarity of the behavior of complex and simple models may in some cases be a coincidence, in other cases it can be formally linked through a rigorous reduction procedure. For example, one of the key problems in climate science is the determination of Equilibrium Climate Sensitivity (ECS) \cite{Charney.1979}. This is defined to be the equilibrium increase in global mean surface temperature that results from a doubling of CO${}_2$ in the atmosphere and as such is the key scientific link between anthropogenic emissions and the climate \cite{knutti2017beyond,ghil2016mathematical}. A rigorous justification of ECS  would need to span the model hierarchy as it involves reduction of ESM runs to an equation that just considers global radiation balance \cite{Hasselmann_1976}.

The aim of this paper to give an overview of prominent and rigorously justifiable reduction methods that have found use in climate modelling, and to highlight others that may find use in future. To keep the review sufficiently brief, we focus on reduction methods for deterministic differential equations and only comment on reduction methods for stochastic differential equations in the outlook. Methods for model reduction have become so large and diverse in recent years that getting an overview of the main approaches for PDEs is already a major challenge when entering the field. Here, we hope to alleviate this issue by providing a relatively concise, yet sufficiently broad, entry point to the recent literature. 

\subsection{Reduction of deterministic climate models}


We present a survey of several major reduction methods that can be applied to climate models. For each approach, we start with its motivation, then present the main reduction idea, giving some concrete examples and comments on the method. In the absence of time-dependent forcing, a climate model can described by
\begin{align}\label{GeneralForm}
 \dot{z}=F(z)
\end{align}
where $z$ represents the discretized state. However, there are slowly evolving aspects (climate/ocean) and faster (weather) variables. This splitting into climate and weather variables, proposed by Hasselmann in \cite{Hasselmann_1976}, has proven to be a very fruitful point of view and have been explored in many papers since; see for example \cite{Majda_Franzke_Crommelin_2009,Franzke_et_al_2015}. Indeed, the wide range of timescales from tectonic to atmospheric processes \cite{Mitchell_1976}, recently reviewed in \cite{Heydt_2021}, highlight there is a multiscale hierarchy inherent in earth system processes. To illustrate this weather-climate scale separation, consider a fast-slow evolution equation of the form 
\begin{align}
 \begin{aligned}\label{HasselmannForm}
 \dot{x}&=f(x,y), \\
 \dot{y}&=\epsilon g(x,y),
 \end{aligned}
\end{align}
where $z=(x,y)$ is split into  slow (climate) variables $y$ and the fast (weather) variables $x$. The separation of timescales is indicated by a small parameter $\epsilon>0$.
Some of the methods, such as slow manifold reductions, averaging or homogenization, try to make use of the special structure such as \eqref{HasselmannForm} with the aim of reducing to the climate variable $y$ only. Others, such as certain Galerkin approximations or asymptotic expansions, can be used to formally establish  \eqref{HasselmannForm} from \eqref{GeneralForm}. 




In Section~\ref{sec:linear}, we consider an abstract forced evolution equation and start with direct linear projection methods, which are among the most classical - yet quite coarse - reduction tools. 
There exist many methods which can help in finding lower-dimensional structures in high dimensional data. In the context of climate data, a reduction via empirical orthogonal functions (EOFs), which we discuss in Section~\ref{Sec:EOFs}, is among the most standard ones.  In Section~\ref{sec:charscales}, we proceed to methods motivated by asymptotic analysis and the principle of dominant balance, where one identifies suitable asymptotic scales directly within evolution equations and formally discards higher-order terms. In Section~\ref{sec:invmfld}, we introduce methods based upon invariant manifold theory, which can be rigorously justified in many cases. 
Note that reduction of the coupled system \eqref{GeneralForm} to a system with the climate variable only does not necessarily lead to a deterministic equation. Sometimes, the variability of the climate system is better described by stochastic forcing or by including a delay term. In fact, the reduction from a deterministic equation to an equation with random forcing can be described by averaging and homogenization techniques which we discuss in Section~\ref{Sec:AvHom}. 
If no direct rigorous deterministic reduction for a given evolution equation is available, one often needs (preliminary) statistical/data-based tools, so we consider classical moment closure schemes in Section~\ref{sec:mc} and indicate with diffusion maps in Section~\ref{Sec:DiffusionMaps} one instance of statistical manifold learning. 


Although we do not claim to cover all possible reduction techniques, we hope that laying out the general ideas via mathematically concrete examples will help the reader better understand any (dis-)advantages for each reduction principle.

\subsection{Classification of Reduction Methods}

Climate models can be described by different types of equations. They can appear as differential equations with possible delay, spatial derivatives and/or stochastic terms. Many of the reduction methods we want to discuss in this article start with one of these types of equations and transform it into another. One can use this to structure the myriad of reduction methods available in the literature. For the few methods we treat in this article, Table~\ref{tab:summary2} give a summary overview of reduction methods explored later in the paper. 

A more precise description of the different methods follow in the subsequent sections, together with some examples of applications to aspects of the climate system.  We hope that this provides readers, who have concrete models in mind, a structured view, which methods are available and how to adapt the main idea to their concrete context. Due to the complexity of climate dynamics and the sort of questions that need to be answered, it is clearly too optimistic to expect a ``one-catch-all method'' for climate model reduction.


\begin{table}
\begin{tabular}{|p{25mm}||c|c|}\hline
    Reduction method & From & To  \\ 
    \hhline{|=#=|=|}
    Galerkin approximations (\S~\ref{Sec:Galerkin}) & PDE & ODE \\ \hline
    \multirow{3}{=}{Empirical orthogonal functions (\S~\ref{Sec:EOFs})} & high dimensional data set & low dimensional data set \\ \cline{2-3}
     & high dimensional ODE & low dimensional ODE \\
     & & 
       \\ \hline
     Principal interaction patterns (\S~\ref{Sec:PIPs}) & high dimensional ODE & low dimensional ODE   \\ \hline
    Characteristic scales (\S~\ref{sec:charscales}) & PDE & PDE of lesser complexity   \\ \hline
    Inertial manifolds (\S~\ref{sec:Inertial}) & PDE & ODE   \\ \hline
    \multirow{3}{=}{Slow manifolds (\S~\ref{sec:Slow})} & \multirow{2}{*}{fast-slow system of PDEs} & PDE with slow variable only   \\ \cline{3-3}
     &  & ODE with slow variable only   \\ \cline{2-3}
     & fast-slow system of ODEs & ODE with slow variable only  \\ \hline
     \multirow{3}{=}{Centre manifolds (\S~\ref{Sec:centre})} & ODE & \multirow{2}{*}{ODE with neutral directions only}   \\ \cline{2-2}
     & PDE &  \\ \cline{2-3}
     & PDE & PDE with neutral directions only   \\ \hline
     \multirow{5}{=}{Averaging (\S~\ref{Sec:Averaging})}  & fast-slow system of ODEs & ODE with slow variable only  \\ \cline{2-3}
      & \multirow{2}{*}{fast-slow system of SDEs} & ODE with slow variable only  \\ \cline{3-3}
      &  & SDE with slow variable only  \\ \cline{2-3}
      & \multirow{2}{*}{fast-slow system of SPDEs} & PDE with slow variable only  \\ \cline{3-3}
      &  & SPDE with slow variable only  \\ \hline
      \multirow{4}{=}{Homogenization (\S~\ref{Sec:Homogenization_Spatially_Periodic}, \ref{Sec:Homog})} & PDE with strongly  & \multirow{2}{*}{PDE with constant coefficients}   \\ 
      & varying coefficients &   \\ \cline{2-3}
      & fast-slow system of ODEs  & SDE with one   \\ 
      & with two time scales & time scale   \\ \hline
      \multirow{2}{=}{Mori-Zwanzig (\S~\ref{Sec:Mori})} & PDE & \multirow{2}{*}{ODE/DDE}   \\ \cline{2-2}
     & ODE &  \\ \hline
     \multirow{2}{=}{Moment closure methods (\S~\ref{sec:mc})} & SDE & \multirow{2}{*}{ODE}   \\ \cline{2-2}
     & PDE &  \\ \hline
     Diffusion maps (\S~\ref{Sec:DiffusionMaps})& high dimensional data set & low dimensional data set   \\ \hline
\end{tabular}
\caption{A summary of the reduction methods discussed in the paper: the relevant Section in the paper is denoted with \S X.Y.}
\label{tab:summary2}
\end{table}





\section{Linear Truncation Methods}
\label{sec:linear}

If one seeks to reduce the complexity of a model or the size of a data set, it is a common approach to project it to a lower-dimensional linear subspace of the state space, which is believed to contain the most important interactions of the model or the data set. The hope is that the parts which are lost during this truncation procedure are not essential and that they can be neglected without losing any meaningful information. Such techniques can be very useful in many different situations: Galerkin approximations can be used to derive well-posedness of partial differential equations on an abstract level~\cite{Evans}, finite element methods are a standard tool in numerical applications~\cite{BrennerScott}, and a principal component analysis is a well established procedure to structure a given data set~\cite{JollifeCadima}. There are many different ways to construct the linear subspace on which the model or data set should be projected and it would be beyond the scope of this paper to discuss all of them. Instead, we just treat three of these methods: First, we study a Fourier-Galerkin approach which can be used to reduce a partial differential equation to an ordinary differential equation. It can also be used to establish a splitting like \eqref{HasselmannForm} in a fast and a slow variable. After that, we discuss empirical orthogonal functions (EOFs) which is a standard technique to reduce the dimensionality of a data set. However, when applied to a dynamical system, it usually only preserves the statistical and not the dynamical properties of the system. There are methods which try to address this issue, one of them being the use of principal interaction patterns (PIPs), which are the last method we study in this section.

\subsection{A Galerkin approach}\label{Sec:Galerkin}

Galerkin approximation appears in different situations in the context of partial differential equations. It is a versatile tool and can be used for deriving the abstract well-posedness of partial differential equations as well as for the development of numerical algorithms. Galerkin approximations are common, and there are many textbooks available, see for example \cite{Fletcher_1984,Thomee_2006}.

\smallskip

{\bf{Idea:}} In certain situations infinite-dimensional systems may come with a natural (Schauder) basis which one can use to construct finite dimensional subspaces. For example, if one has a partial differential equation with periodic boundary conditions, this can be the canonical Fourier basis, or if there is a self-adjoint operator on a Hilbert space involved, then this may be the set of eigenfunctions of this operator. This is a Galerkin kind of approach, which is often most convenient if periodic boundary conditions in space are imposed. The idea is to expand the equation into a Fourier series in space and to omit high frequencies. Oftentimes, such a basis comes with a natural order. In the case of the Fourier basis, this order is given by the increasing frequencies, and in the case of eigenfunctions of a self-adjoint operator, it may be given by the order of the eigenvalues on the real axis. The finite-dimensional approximations are then obtained by projecting the system to the linear span of a finite number of basis vectors.

\smallskip

{\bf{Outline of the Reduction Procedure:}} A standard scenario in which Galerkin approximations can be applied is given by the abstract evolution equation
\begin{align}\label{Eq:Galerkin_EvolutionEq}
\partial_t u(t) = F(t,u(t)),\quad u(0)=u_0
\end{align}
in a Hilbert space $H$. Let $(e_k)_{k\in\N}$ be an orthonormal basis of $H$. Then one can define the $n$-dimensional spaces $H_n:=\operatorname{span}\{e_1,\ldots, e_n\}$ $(n\in\N)$ and the canonical orthogonal projections $P_n\colon H\to H_n$. If we truncate \eqref{Eq:Galerkin_EvolutionEq} to $H_n$, we obtain the finite-dimensional system
\begin{align}\label{Eq:Galerkin_EvolutionEq_Projected}
  \partial_t v_n(t) = P_n F(t,v_n(t)),\quad v_n(0)=P_n u_0,
\end{align}
where $v_n:=P_n u$. Let $\langle\,\cdot,\cdot\rangle$ be the inner product of $H$. Then we can rewrite \eqref{Eq:Galerkin_EvolutionEq_Projected} in terms of
\begin{align}\label{Eq:Galerkin_EvolutionEq_Projected2}
 \partial_t  \sum_{j=1}^n\langle v_n(t),e_j\rangle e_j = \sum_{j=1}^n \langle P_n F(t,v_n(t)),e_j\rangle e_j,\quad \sum_{j=1}^n\langle v_n(t),e_j\rangle e_j =\sum_{j=1}^n\langle u_0,e_j\rangle e_j.
\end{align}
Now, we define $x_{j}^n(t):=\langle v_n(t),e_j\rangle$, $f^n_j(t,x_1^n(t),\ldots,x_n^n(t)):=\langle P_n F(t,v_n(t)),e_j\rangle$ and $x_{0,j}^n:=\langle u_0,e_j\rangle$ so that \eqref{Eq:Galerkin_EvolutionEq_Projected2} can be rewritten as
\begin{align}
\frac{d}{dt}x_j^n (t)= f^n_j(t,x_1^n(t),\ldots,x_n^n(t)),\quad x_j^n (0)=x_{0,j}^n,\quad (j=1,\ldots,n)
\end{align}
by the linear independence of $(e_k)_{k\in\N}$. If we define 
\begin{align*}
    x^n&:=(x_1^n,\ldots,x_n^n),\\
    f^n&:= (f_1^n,\ldots,f_n^n),\\
    x_0^n&:=(x_{0,1}^n,\ldots,x_{0,n}^n),
\end{align*} 
then we obtain the more compact form
\begin{align}
\dot{x}^n(t) = f^n(t,x^n(t)),\quad x^n (0)=x_{0}^n.
\end{align}
This is now an $n$-dimensional system of ordinary differential equations and it is called Galerkin approximation of \eqref{Eq:Galerkin_EvolutionEq}. Moreover, under relatively mild conditions its holds that $x^n\to x$ weakly as $n\to\infty$ in $H$. For more details on the precise assumptions and statements we refer for example to \cite[Theorem 2.76]{Ruzicka_2020} or \cite[Chapter III.4]{Showalter_1997}.
  

\smallskip

 {\bf{Example:}} We follow \cite{Majda_Timofeyev_VandenEijnden_2001}, where the equations for barotropic flow on a beta plane with topography and mean flow are considered:
 \begin{equation}
    \begin{aligned}\label{Eq:BarotropicFlow1}
    \frac{\partial q}{\partial t}+\nabla^{\bot}\psi\cdot\nabla q + U\frac{\partial q}{\partial x}+\beta\frac{\partial\psi}{\partial x}&=0,\\
    q&=\Delta\psi+h,\\
    \frac{dU}{dt}&=\dashint h\frac{\partial\psi}{\partial x}.
    \end{aligned}
 \end{equation}
 Here, $q(x,y,t)$ denoted the small-scale potential vorticity, $U(t)$ is the mean flow, $\psi(x,y,t)$ is the small-scale stream function, and $h(x,y)$ denotes the underlying topography, whereas $\beta$ approximates the variation of the Coriolis parameter. In \cite{Majda_Timofeyev_VandenEijnden_2001}, the equation is studied with periodic boundary conditions so that one may take the Fourier basis functions $([x\mapsto e^{ikx}])_{k\in\Z^2}$ as an orthonormal basis of $L_2(\mathbb{T}^2)$. As in \cite[Section 3]{Majda_Timofeyev_VandenEijnden_2001} we choose $\Lambda\in\N$ and define $\sigma_{\Lambda}=\{k\in\Z^2: 1\leq|k|\leq\Lambda\}$, $B_\Lambda=\{e^{ikx}:k\in\sigma_{\Lambda}\}$ as well as the orthogonal projection $P_{\Lambda}$ onto $H_{\Lambda}:=\operatorname{span} B_{\Lambda}$. Moreover, we write
 \begin{align*}
  \psi_{\Lambda}(x)&=\sum_{1\leq|k|^2\leq\Lambda} \hat{\psi}_k(t)e^{ikx},\quad h_{\Lambda}(x)=\sum_{1\leq|k|^2\leq\Lambda} \hat{h}_ke^{ikx},\\ q_{\Lambda}(x,t)&=\sum_{1\leq|k|^2\leq\Lambda}\hat{q}_k(t)e^{ikx} 
 \end{align*}
for the truncations to $H_{\Lambda}$ and substitute $u_k(t):=|k|\psi_k(t)$. By the reduction procedure described above, one obtains the system
\begin{align*}
     \frac{dU}{dt} &= \operatorname{Im} \sum_{k\in\sigma_{\Lambda}} H_k \overline{u_k},\\
     \hat{q}_k&=-|k|u_k+\hat{h}_k \quad(k\in\sigma_{\Lambda}),\\
     \frac{du_k}{dt}&=iH_kU-i(k_xU-\Omega_k)u_k+\sum_{l\in\sigma_{\Lambda}}L_{kl}u_l+\frac{1}{2}\sum_{
    \genfrac{}{}{0pt}{}{l,m
     \in\sigma_{\Lambda}}{k+l+m=0}} B_{klm}
     \overline{u}_l \overline{u}_m\quad(k\in\sigma_{\Lambda}),
\end{align*}
where
\begin{align*}
    L_{kl}&=\frac{k_xl_y-k_yl_x}{|k||l|}h_{k-l},\quad B_{klm} = (l_ym_x-l_xm_y)\frac{|l|^2-|m|^2}{|k||l||m|},\\   \Omega_k&=\frac{k_x\beta}{|k|^2},\quad H_k=\frac{k_x\hat{h}_k}{|k|}.
\end{align*}
The Fourier basis could also be used to derive a splitting in a weather and a climate variable as in \eqref{HasselmannForm}, since higher Fourier modes have stronger oscillations which oftentimes lead to faster dynamics. In this sense, one can introduce $1\leq \bar{\Lambda}\leq\Lambda$ and consider $u_k,q_k$ with $k\in\sigma_{\overline{\Lambda}}$ and $U$ as slow climate variables and  $u_k,q_k$ with $k\in\sigma_{\Lambda}\setminus\sigma_{\overline{\Lambda}}$ as fast weather variables. We refer to \cite{Majda_Timofeyev_VandenEijnden_2001} for further details and to Section \ref{Sec:AvHom} where we continue with this example.
  
\smallskip

 {\bf{Comments:}} 
 
 \begin{itemize}
  \item Galerkin-type methods have the advantage that they are often easy to implement and that they provide a nice geometric intuition as solutions are expanded in basis functions. Therefore, one hopes to capture most features that can be expressed in terms of the finite number of basis functions left in the approximation. Numerical tests can often indicate, how many basis functions are needed in practice.
  \item One disadvantage is that models can become large very quickly if a sufficient accuracy is required, which somewhat defies the initial point of model reduction as a system of nonlinear ordinary differential equation with tens or hundreds of variables is analytically frequently not much easier than directly analysing the initial partial differential equation. Particularly for fluid dynamics, e.g., Navier-Stokes-type equations, it is well-known that small-scale turbulence can spread from very high Fourier modes into low modes in certain cases.
\end{itemize}

In fact, we shall see below that the advantages and disadvantages mentioned above, are quite typical not only for Galerkin-type methods but for all methods based upon linear projection/reduction. The next method shows another classical example of this principle.

\subsection{Empirical orthogonal functions (EOFs)}
\label{Sec:EOFs}

There are many different names for the method that yields the construction of the so-called empirical orthogonal functions. Closely related methods are also referred to as singular value decomposition \cite{bretherton1992svd}, principal component analysis \cite{preisendorfer1988pca} or proper orthogonal decomposition \cite{chen2019pod}. This technique is well established and widely used, not only for the reduction of climate models or climate data. We refer to \cite{Bjornsson_1997,Hannachi_et_al_2007,Lorenz_1956,Navarra_Simoncini_2010}.

\smallskip

 {\bf{Idea:}} The aim is to construct a subspace of a given dimension which carries the largest amount of variance of the given data set. The data set is then projected orthogonally to this subspace. If one starts with a dynamical system instead of a data set, then one can obtain a data set by taking a certain number of snapshots of trajectories of the dynamical system.
 
\smallskip

 {\bf{Outline of the Reduction Procedure:}} One starts with a collection of $N\in\N$ points $x_1,\ldots, x_N\in\R^d$. Alternatively, if a dynamical system in a $d$-dimensional state space is given, then one can take $N$ snapshots of a trajectory $(x_1,\ldots, x_n)=(z(t_1),\ldots, z(t_N))$ at different times $t_1,\ldots, t_N\geq0$ or a collection of such snapshots of trajectories with different initial values. Either way, one works with a matrix $X=(x_1-\bar{x},\ldots, x_N-\bar{x})\in \R^{d\times N}$ which is given as a concatenation of the data set relative to the mean $\bar{x}:=\frac{1}{N}\sum_{j=1}^N x_j$. The matrix $XX^T\in\R^{d\times d}$ is symmetric and positive semidefinite. Therefore, it has $d$ eigenvalues $\lambda_1\geq\ldots \geq \lambda_d\geq0$ and corresponding orthonormal eigenvectors $v_1,\ldots v_d\in\R^d$. From a statistical point of view the matrix is the covariance matrix of the data set multiplied $N$ (or $N-1$ depending on the convention). Thus, the eigenvalue $\lambda_i$ can be considered as a measure of variance the direction of $v_i$ carries. Since the initial idea was to determine the directions with the most variance, one projects the data set to the subspace generated by $v_1,\ldots,v_k$ for some $k\leq d$. More precisely, one defines
 \[
	x^k_l:= \sum_{j=1}^k \langle x_l, v_j \rangle v_j\quad (l=1,\ldots, N).
 \]
 Here, $\langle x_l, v_j \rangle$ denotes the usual Euclidean scalar product, but note that also other choices of scalar products are possible, as we will explain later in the comments. The points $x^k_1,\ldots, x^k_N$ are now contained in a $k$-dimensional subspace of $\R^d$. The choice of $k$ depends on the application in mind and can also be chosen according to the given data set. Small values of $k$ have the advantage that one ends up with a very low-dimensional data set. However, the approximation gets rougher with smaller values of $k$. Hence, there is a trade-off between the dimension of the reduced data set and the quality of approximation of the original data set.
  
\smallskip

 {\bf{Example:}} There are many works that use EOFs and its extensions in a climate context, see for example \cite{Beckers_Rixen_2003,Biau_et_al_1999,Crommelin_Majda_2004,Kawamura_1994,Kessler_2001,Selten_1997,Smith_et_al_1996,Wallace_et_al_1993}. Let us briefly explain the example in \cite{Crommelin_Majda_2004}. Therein, the authors compare reductions of the Charney--DeVore model with different basis functions, among them EOFs. The Charney--DeVore model was derived in \cite{DeSwart_1989} as a $6$-dimensional Galerkin approximation of the vorticity equation for a large-scale atmospheric flow on a $\beta$-plane channel with topography which similar to equation \eqref{Eq:BarotropicFlow1}. The reduced equations which are considered in \cite{Crommelin_Majda_2004} are of the form
 \begin{align}\label{Eq:CDV_model}
     \begin{aligned}
        \dot{x}_1 &= \tilde{\gamma}_1x_3 - C (x_1-x_1^{\ast}),\\
        \dot{x}_2 &= -(\alpha_1x_1-\beta_1)x_3-Cx_2-\delta_1x_4x_6,\\
        \dot{x}_3 &= (\alpha_1x_1-\beta_1)x_2-\gamma_1x_1-Cx_3+\delta_1 x_4x_5,\\
        \dot{x}_4 &= \tilde{\gamma}_2x_6-C(x_4-x_4^{\ast})+\eta(x_2x_6-x_3x_5),\\
        \dot{x}_5 &= -(\alpha_2x_1-\beta_2)x_6-Cx_5-\delta_2x_4x_3,\\
        \dot{x}_6 &= (\alpha_2x_1-\beta_2)x_5-\gamma_2x_4-Cx_6+\delta_2x_4x_2,\\
     \end{aligned}
 \end{align}
 where $(x_1,x_2,x_3,x_4,x_5,x_6)$ are the unknown functions, $x_1^{\ast},x_4^{\ast}$ are forcing terms and the other quantities are model parameters taking the form
 \begin{align*}
    &\alpha_j=\frac{8\sqrt{2}}{\pi}\frac{j^2}{4j^2-1}\frac{b^2+j^2-1}{b^2+j^2},\quad\beta_j=\frac{\beta b^2}{b^2+j^2},\quad \delta_j=\frac{64\sqrt{2}}{15\pi}\frac{b^2-j^2+1}{b^2+j^2},\\
    &\tilde{\gamma}_j=\gamma\frac{4j}{4j^2-1}\frac{\sqrt{2}b}{\pi},\quad\eta=\frac{16\sqrt{2}}{5\pi},\quad\gamma_j=\gamma\frac{4j^3}{4j^2-1}\frac{\sqrt{2}b}{\pi(b^2+j^2)}\quad(j\in\{1,2\}).
 \end{align*}
 Here, $b$ models the length-width ration of the $\beta$-channel, $\beta$ comes from the Coriolis force, $\gamma$ from the topographic height and the damping parameter $C$ is determined by the friction in the Ekman layer. In the vorticity equation from which \eqref{Eq:CDV_model} is derived, it is assumed that the forcing term points in zonal direction only. During the Galerkin reduction, this leads to the fact that there are only forcing terms in the equations for $x_1$ and $x_4$. This model can show rapid and chaotic transitions between different flow regimes. As was shown in \cite{Crommelin_Opsteegh_Verhulst_2004}, this is the case for the parameter values $(x_1^{\ast},x_4^{\ast},C,\beta,\gamma,b)=(0.95,\;-0.76095,\;0.1,\;1.25,\;0.2,\;0.5)$. In \cite{Crommelin_Majda_2004} it was investigated, whether these regime transitions would persist under an EOF reduction. As described in the outline of the EOF reduction procedure above, the EOFs were computed from the data set generated by a numerical integration of the full model. The authors also used different energy metrics, i.e., different scalar products for the computation of the EOFs. It was found that none of the reduced models was able to reproduce the chaotic transition behavior of the full model. Even though some of the reduced models show regime transitions, they are far too regular and of periodic type. This already indicates one of the problems of EOF reductions. Even though reduced models usually shows the same or a similar statistical behavior as the original model, the dynamical behavior might still be very different.\smallskip

 {\bf{Comments:}}  \begin{itemize}
  \item In applications, variance and distances are sometimes not measured in terms of the Euclidean metric, but with respect to a certain energy metric which is given as a symmetric positive definite matrix $M\in\R^{d\times d}$. In this case, orthogonality is not understood with respect to the Euclidean scalar product, but with respect to the scalar product
  \[
  	(x,y)_M=\langle x, My\rangle.
  \]
  In this case, the covariance matrix is given by $M^{1/2}X(M^{1/2}X)^\top$ and $\lambda_1^{(M)}$, $\ldots$, $\lambda_N^{(M)}$ as well as $\tilde{v}_1^{(M)},\ldots, \tilde{v}_N^{(M)}$ are defined by the relation
  \[
  	M^{1/2}X(M^{1/2}X)^T\tilde{v}_i^{(M)}=\lambda_i^{(M)}\tilde{v}_i^{(M)}.
  \]
  As above, the $\tilde{v}_1^{(M)},\ldots, \tilde{v}_N^{(M)}$ are orthonormal with respect to the Euclidean scalar product. Therefore, the basis defined by $v_i^{(M)}:=M^{-\frac{1}{2}}\tilde{v}_N^{(M)}$ is orthonormal with respect to $M$, i.e.
  \[
  	(v_i^{(M)},v_j^{(M)})_M:=\langle v_i^{(M)}, Mv_j^{(M)}\rangle=\delta_{i,j}.
  \]
  The projection of the data points is then given by 
  \[
  	x^{k,(M)}_l:= \sum_{j=1}^k ( x_l, v_j^{(M)} )_M v_j^{(M)}=\sum_{j=1}^k \langle x_l,M v_j^{(M)} \rangle v_j^{(M)}
  \]
  where $l=1,\ldots, N$.
  In the literature, it is often common to define $v_1^{(M)},\ldots, v_N^{(M)}$ directly as the eigenvectors of the eigenvalue problem
  \[
  	XX^\top Mv_i^{(M)}=\lambda_i^{(M)}v_i^{(M)},
  \]
  see for example \cite[Section 3]{Crommelin_Majda_2004} or \cite[Section 4.3]{Kwasniok_1996}. But it is easy to see that both definitions are mathematically equivalent. Yet, it is often natural to re-scale via an energy metric to bring different directions onto a common scale. In fact, introducing the matrix $M$ into the inner product corresponds to an an-isotropic or ellipsoidal deformation of phase space. This is a common trick not only appearing from a physical or geometric perspective but also in the context of statistical techniques, where the Mahalanobis distance uses a scaling matrix $M$ based upon the covariance of a stochastic process~\cite{Mahalanobis,KuehnSDEcont1}. 
  \item In some sense adding the energy metric $M$ is an extension of the usual EOF approach. However, this is not the only possible extension. We refer the reader to \cite{Preisendorfer_1988} and \cite{Hannachi_et_al_2007} for more content in this direction.
  \item EOFs are designed to reduce the dimension of a data set while giving a good approximation of its statistical features. However, this method does not aim to preserve the dynamics of a dynamical system and indeed it has problems in doing so. This is pointed out in several papers like \cite{Crommelin_Majda_2004,Monahan_et_al_2009} and references therein. In particular, one can not expect that dynamical notions such as invariant manifolds, attractors or bifurcations persist under an EOF-reduction.
 \end{itemize}

 As discussed above, EOFs can not be expected to preserve the dynamics of a dynamical system in general. A reduction method via principal interaction patterns, which we discuss next, is designed to improve on this issue.

\subsection{Principal interaction patterns (PIPs)}
\label{Sec:PIPs}

Now we consider another choice of basis functions, namely principal interaction patterns (PIPs). The idea behind PIPs was introduced by Hasselmann in \cite{Hasselmann_1988}. It was then used and further refined in several papers, e.g. \cite{Achatz_Schmitz_1997,Achatz_Schmitz_Greisinger_1995, Kwasniok_1996,Kwasniok_1997,Kwasniok_2001,Kwasniok_2004}.

\smallskip

 {\bf{Idea:}} The idea is similar to the one of the EOF-method in the sense that one is looking for a suitable low-dimensional subspace to which one can project a dynamical system. But in contrast to the EOF-method, principal interaction patterns were introduced with the aim of finding a reduction method which focuses more on the dynamics instead of just on the statistics of a dynamical system. Therefore, PIPs are constructed as (approximate) minimizers of a certain error functional which measures the distance of the projected solutions to the original solutions.
 
\smallskip

 {\bf{Outline of the Reduction Procedure:}}  The aim of this method is to reduce a system of ordinary differential equations
 \begin{align}\label{Eq:PIP_ODE}
 	\dot{x}(t)=F(t,x(t)),\quad x(0)=x^0
 \end{align}
 in the state space $\R^d$ with $d\in\N$ being large to a lower-dimensional system. Even though PIPs are constructed for systems as \eqref{Eq:PIP_ODE}, they can also be used for partial differential equations of the form
 \begin{align}\label{Eq:PIP_PDE}
  \partial_tu=D(u)
 \end{align}
 with a suitable nonlinear operator $D$. But in such a case, \eqref{Eq:PIP_PDE} is approximated by \eqref{Eq:PIP_ODE} through a Galerkin approximation and then the PIPs are constructed for \eqref{Eq:PIP_ODE}.\\
 For the construction of PIPs one needs a suitable notion of an average or expectation of several quantities over a set of inital values \eqref{Eq:PIP_ODE}. We will use the notation $\langle \cdot \rangle$ for this average. In the literature, this average is not always made precise and may slightly differ in different situations. The idea behind the precise definitions is however similar in most cases. In the original paper by Hasselmann \cite{Hasselmann_1988}, in which the idea of PIPs was first introduced, it is just described as expectation value. In \cite{Kwasniok_1997} it is an ensemble average over an attractor. Similarly, an ensemble average over a finite number of such initial conditions was chosen in \cite{Crommelin_Majda_2004}. Since PIPs are eventually determined numerically anyway, we also work with a finite number of initial conditions in the following.
Therefore, let $I=\{x^{0,1},\ldots, x^{0,N}\}$ be a set of $N\in\N$ initial conditions for \eqref{Eq:PIP_ODE}. Whenever $V$ is a vector space and $f\colon\R^d\to V$ a mapping, we use the notation
\[
	\langle f \rangle := \frac{1}{N}\sum_{j=1}^N f(x^{0,j}).
\]
But one should keep in mind that other notions of average are also possible.\\
For a fixed number $k\in \N$ with $k< d$ and linearly independent vectors $p_1,\ldots, p_k\in\R^d$ one defines $P:=\operatorname{span}\{p_1,\ldots,p_k\}$. Then, one projects \eqref{Eq:PIP_ODE} orthogonally to the reduced system
 \begin{align}\label{Eq:PIP_ODE_projected}
	\dot{x}_P(t)=F_P(t,x_P(t)),\quad x_P(0)=x^0_P.
 \end{align}
 Here, $F_P=\operatorname{pr}_P \circ F$ where $\operatorname{pr}_P$ denotes the orthogonal projection onto $P$ and where orthogonality may again be understood with respect to the scalar product $\langle x, My\rangle$ given by the Euclidean scalar product and the energy metric $M\in\R^{d\times d}$.\\
 Now, we split 
 	\[
 		x(t)=\langle x(t)\rangle+\tilde{x}(t),\quad x_P(t)= \langle x(t)\rangle+\tilde{x}_P(t)
 	\]
 	and work with the deviations from the mean state $\langle x(t)\rangle$ of the system at time $t$. We define the error functional
  \[
  Q(P;w):=\left\langle\int_{0}^{\infty} \| \tilde{x}_P(\tau)-\tilde{x}(\tau)\|_{M}^2 w(\tau)\,d\tau\right\rangle.
 \]
 Here, $\|v\|_M:= \sqrt{\langle v,Mv\rangle}$ and $w\colon[0,\infty)\to [0,\infty)$ is a suitable weight function. As for the definition of the average, there is some freedom in the choice of the error functional that goes beyond different choices of $w$ or $M$. We refer the reader to \cite{Crommelin_Majda_2004,Hasselmann_1988,Kwasniok_1997,Kwasniok_2004}, where different choices have been made. Now, $P$ is chosen as an approximate minimizer of the error functional $Q(P;w)$ and usually determined by numerical methods. Details of the minimization procedure are carried out in the appendix of \cite{Kwasniok_1997}. The principal interaction patterns are now the basis $p_1,\ldots,p_k$ of $P$. Obviously, they are not unique, since one can obtain a new basis by a linear transformation on $P$. Thus, one usually imposes constraints on the choice of $p_1,\ldots,p_k$ so that they can be uniquely determined, see for example \cite[Section II.C]{Kwasniok_1997}.
  
\smallskip

 {\bf{Example:}} We consider again the CDV-model \eqref{Eq:CDV_model}. In the section on EOF reductions we discussed the reductions carried out in \cite{Crommelin_Majda_2004} and how an EOF reduction could not reproduce the chaotic transitions between different flow regimes. Since PIPs were introdruced with the intention in mind, that they should be able to preserve the dynamical behavior of a system, one might ask whether they are able to do so in the case of the CDV-model. This was also investigated in \cite{Crommelin_Majda_2004}. The authors carried out various PIP reductions with different weights $w$ and different sets of initial conditions $I$. It was found that most models were able to reproduce the chaotic regime transitions. Moreover, the power spectra of the trajectories integrated from the reduced models were compared to the ones from the full model and it was found that they look very similar. Therefore, one could say that for the CDV-model a PIP reduction would indeed yield a more precise reduced model than an EOF reduction. However, the dependence of the reduction on the input parameters is in some cases counterintuitive. In particular, it was found that the chaotic regime transitions would not persist if the set of initial values $I$ was too large. This shows that even if a PIP reduction can produce very good results, but the results depend crucially on the input parameters and it is not always clear, what choice of parameters would be suitable. \smallskip

 {\bf{Comments:}} \begin{itemize}
                 \item The choice of $w$, $I$ and the averaging procedure is crucial for the outcome of the minimization problem. This can be seen as an advantage of PIPs, as it provides some flexibility. But the downside is that it has not yet been systematically investigated which choices would yield good results for a given model. The outcome can also be counterintuitive to some extend. For example, as discussed above it was observed  in \cite{Crommelin_Majda_2004} that for the Charney-DeVore model PIPs have difficulties in preserving the dynamics if the set $I$ contains too many initial values.
                 \item For example \cite[Section II.E]{Kwasniok_1997} it was observed that one can obtain EOFs as a limiting case of PIPs if $w=\frac{1}{\tau_{\max}}\one_{\tau_{\max}}$ and $\tau_{\max}\to0$.
                 \item The minimization of $Q(P;w)$ resembles problems encountered in optimal control theory as for example presented in \cite{Lions_1971}. The problem here is slightly different in the sense that one is looking for a subspace such that the projection error is minimized in terms of $Q$, whereas in optimal control one usually looks for a control function or parameter such that a cost functional is minimized subject to some constraints. But nonetheless, techniques from optimal control also help in numerically solving the minimization task for $Q$, see \cite[Appendix B]{Kwasniok_1997}.
                \end{itemize}

In the previous section we have seen how partial differential equations can by reduced to finite-dimensional ordinary differential equations by projecting them to finite-dimensional subspaces of the phase space. The knowledge of certain patterns or characteristic parameters may be useful for such an approach, but it is not essential. This is different in the next method we are going the describe, as it is based on the utilization of the characteristic length scales, time scales and amplitudes of the underlying system.

\section{Characteristic Scales}
\label{sec:charscales}

Characteristic length or time scales can sometimes be found in evolution equations of the form
\begin{align}\label{Eq:CharScales}
	\partial_t u(t,x) = F(\partial_x,u(t,x),\epsilon),
\end{align}
where $F$ denotes a generally nonlinear spatial differential operator and $\epsilon$ denotes a small parameter. The notion of characteristic scale is not defined in a rigorous mathematical manner but is instead based on certain parameters or observations of the system described by \eqref{Eq:CharScales}.
 
\smallskip

 {\bf{Idea:}}  Even very complex dynamical systems may exhibit characteristic length or time scales at which interesting patterns emerge. Such characteristic scales can for example be induced by certain parameters in the system or sometimes they are just observed as the system evolves. If one knows the underlying equations of the system, then one can rescale them according to the characteristic scales and just consider the leading order terms which appear after the transformation. The hope is that by restricting to the leading order terms the equations become simpler to analyze while their solutions might still show the same patterns whose scales were considered as characteristic. Effectively, this idea is directly motivated by classical formal matching techniques for differential equations~\cite{BenderOrszag,DeJagerFuru,KevorkianCole}, where one also searches for scalings such that several terms within an equation are in a dominant balance, allowing us to discard higher-order terms based on asymptotic analysis.
  
\smallskip

 {\bf{Outline of the Reduction Procedure:}}   We suppose that \eqref{Eq:CharScales} is given in in a reference time and reference length scale in which we denote the time variable by $t$ the space variable by $x$. Moreover, there are $M\in\N$ characteristic time and $N\in\N$ characteristic length scales given by the coordinate transformations $\tau_{(j)}= \epsilon^{\alpha_{t,j}} t$, $(j=1,\ldots,M)$, and $\xi_{(k)}=\epsilon^{\alpha_{x,k}}x$ $(k=1,\ldots,N)$, where the $\alpha_{t,j}$ and $\alpha_{x,j}$ are given real numbers. Let us fix $j_0\in\{1,\ldots,M\}$ and $k_0\in\{1,\ldots,N\}$ and suppose that we are interested in patterns that emerge on the scales given by the variables $(\tau_{(j_0)},\xi_{(k_0)})$. Therefore, we rescale \eqref{Eq:CharScales} according to this change of variables and obtain
 \begin{align}\label{Eq:CharScales:Rescaled}
	 \epsilon^{\alpha_{t,j_0}}\partial_{\tau_{(j)}} u(\tau_{(j)},\xi_{(k_0)})=F(\epsilon^{\alpha_{x,k_0}}\partial_{\xi_{(k_0)}},u(\tau_{(j)},\xi_{(k_0)}),\epsilon).
 \end{align}
 Depending on the exponents $\alpha_{t,j_0}$ and $\epsilon^{\alpha_{x,k_0}}$ one has different leading order terms, i.e., terms with the lowest power in $\epsilon$. At this point, it can already be insightful to study the equation one obtains by omitting the higher order terms. But in practice one sometimes also has more information about the amplitudes of certain patterns so that one can write
 \[
 	u(\tau_{(j)},\xi_{(k_0)})=\sum_{l=1}^P \epsilon^{\beta_l} u_l(\tau_{(j)},\xi_{(k_0)})
 \]
 for certain $P\in\N$, $\beta_1,\ldots,\beta_P\in\R$ and potentially mutually dependent functions $u_1,\ldots,u_P$. Substituting this expansion into \eqref{Eq:CharScales:Rescaled} may yield different leading order terms and thus also different reduced equations.  Mathematically, including amplitudes is directly motivated by the theory of amplitude/modulation equations~\cite{CrossHohenberg,Hoyle,KuehnBook1} used near bifurcations. For ODEs one can hope to making the scaling arguments geometrically rigorous, even near bifurcations, via geometric desingularization~\cite{DumortierRoussarie,JardonKuehn,KruSzm3,Kuehn_2015}. 
   
\smallskip

 {\bf{Example:}}  In \cite{Klein_2010}, formal asymptotic analysis ideas regarding characteristic are applied to compressible flow equations for an ideal gas with constant
specific-heat capacities and including gravity, rotation, and generalized source terms. In addition, a tangent plane approximation is used, i.e. one negelects the curvature of the sphere but still takes the Coriolis force into account. More precisely, the system
\begin{align}\label{Eq:CharScales:Example}
	\begin{aligned}
 \left(\frac{\partial}{\partial t}+v_{\parallel}\cdot\nabla_{\parallel}+w\frac{\partial}{\partial z}\right)v_{\parallel}+\epsilon (2\Omega\times v)_{\parallel}+\frac{1}{\epsilon^3\rho}\nabla_{\parallel}p&=Q_{v_\parallel},\\
  \left(\frac{\partial}{\partial t}+v_{\parallel}\cdot\nabla_{\parallel}+w\frac{\partial}{\partial z}\right)w+\epsilon(2\Omega\times v)_{\bot}+\frac{1}{\epsilon^3\rho}\frac{\partial p}{\partial z}&=Q_w-\frac{1}{\epsilon^3},\\
  \left(\frac{\partial}{\partial t}+v_{\parallel}\cdot\nabla_{\parallel}+w\frac{\partial}{\partial z}\right)\rho+\rho\nabla\cdot v&=0,\\
  \left(\frac{\partial}{\partial t}+v_{\parallel}\cdot\nabla_{\parallel}+w\frac{\partial}{\partial z}\right)\Theta&=Q_{\Theta}
  \end{aligned}
\end{align}
is considered. In this system, $v_{\parallel}$ and $w$ denote (dimensionless) horizontal and vertical flow velocity, p the pressure, $\rho$ the density, $\Theta=p^{1/\gamma}/\rho$ the potential temperature (with $\gamma$ being the ``dry isentropic exponent''), $\Omega$ the earth rotation vector and $Q_{[\cdot]}$ some source terms. Moreover, $z$ corresponds to the vertical direction, $x$ to the horizontal direction and $t$ to the time. \\
  First, the pressure $p$ is replaced by the ``Exner pressure'' $\pi$, i.e. $\pi=p^{\Gamma}$ with $\Gamma=\frac{\gamma-1}{\gamma}$. Now, time and horizontal coordinates are rescaled by $\tau=\epsilon^{\alpha_t}x$ and $\xi=\epsilon^{\alpha_x}x$ with certain $\alpha_t,\alpha_x\geq0$ which depend on the scale one wants to consider. Fluctuations of the potential temperature are expected to be of order $\epsilon$, while velocities in horizontal directions scale as $\epsilon^{\alpha_x}$. Therefore, one may introduce the expansions
 \begin{align*}
  \Theta(\tau,\xi,z)&=1+\epsilon \bar{\theta}(z)+\epsilon^{\alpha_{\pi}}\tilde{\theta}(\tau,\xi,z),\\
  \pi(\tau,\xi,z)&=\bar{\pi}(z)+\epsilon^{\alpha_\pi}\Gamma\tilde{\pi}(\tau,\xi,z),\\
  w(\tau,\xi,z)&=\epsilon^{\alpha_x}\tilde{w}(\tau,\xi,z).
 \end{align*}
Here, $\bar{\pi}$ satisfies
\[
 1-\Gamma\int_0^z\frac{1}{\bar{\Theta}(z')}\,dz'
\]
with $\bar{\Theta}=1+\bar{\theta}$ being the horizontally averaged mean stratification of the potential temperature. Moreover, the parameter $\alpha_{\pi}$ has to be chosen later to balance the pressure gradient in the horizontal momentum equation, i.e., the first equation of \eqref{Eq:CharScales:ExampleRescaled}. Using these substitutions, the system \eqref{Eq:CharScales:Example} can be rewritten as
\begin{align}\label{Eq:CharScales:ExampleRescaled}
\begin{aligned}
 \left(\frac{\epsilon^{\alpha_t}}{\epsilon^{\alpha_x}}\frac{\partial}{\partial \tau}+v_{\parallel}\cdot\nabla_{\xi}+\tilde{w}\frac{\partial}{\partial z}\right)v_{\parallel}+ \frac{(2\Omega\times v)_{\parallel}}{\epsilon^{\alpha_x-1}}+\frac{\epsilon^{\alpha_{\pi}}}{\epsilon^3}\Theta\nabla_{\xi}\tilde{\pi}&=Q_{v_\parallel}^{\epsilon},\\
 \left(\frac{\epsilon^{\alpha_t}}{\epsilon^{\alpha_x}}\frac{\partial}{\partial \tau}+v_{\parallel}\cdot\nabla_{\xi}+\tilde{w}\frac{\partial}{\partial z}\right)\tilde{w}+\frac{(2\Omega\times v)_{\bot}}{\epsilon^{2\alpha_x-1}}+\frac{\epsilon^{\alpha_\pi}}{\epsilon^{3+2\alpha_x}}\left(\Theta\frac{\partial\tilde{\pi}}{\partial z}-\frac{\tilde{\theta}}{\Theta}\right)&=Q_w^{\epsilon},\\
  \left(\frac{\epsilon^{\alpha_t}}{\epsilon^{\alpha_x}}\frac{\partial}{\partial \tau}+v_{\parallel}\cdot\nabla_{\xi}+\tilde{w}\frac{\partial}{\partial z}\right)\tilde{\pi}+\frac{\gamma\pi}{\epsilon^{\alpha_\pi}}\left(\nabla_{\xi}\cdot v_{\parallel}+\frac{\partial\tilde{w}}{\partial z}+\frac{\tilde{w}}{\gamma\Gamma\pi}\frac{d\tilde{\pi}}{dz}\right)&=\frac{\gamma\pi Q^{\epsilon}_{\Theta}}{\Theta},\\
  \left(\frac{\epsilon^{\alpha_t}}{\epsilon^{\alpha_x}}\frac{\partial}{\partial \tau}+v_{\parallel}\cdot\nabla_{\xi}+\tilde{w}\frac{\partial}{\partial z}\right)\tilde{\theta}+\frac{\epsilon}{\epsilon^{\alpha_\pi}}\tilde{w}\frac{d\tilde{\theta}}{dz}&=Q_{\Theta}^{\epsilon}.
  \end{aligned}
\end{align}
For different choices of $\alpha_x,\alpha_t,\alpha_{\pi}$ one obtains different leading order equations. In \cite{Klein_2010}, it is assumed that the order of the source terms is high enough, so that they can be neglected.
   
\smallskip

 {\bf{Comments:}}\begin{itemize}
                 \item In \cite[Figure 1]{Klein_2010} it is shown which length and time scales correspond to which models in the case of atmospheric flows. Roughly speaking, a formal matching of characteristic scales induces a hierarchy of models, where it is crucial to cross-check final analysis results regarding their validity at the next lower or higher scale. 
                 \item The reduction in this section is a formal approach based on asymptotic expansions. Although these expansions can sometimes be made rigorous, it is often difficult to prove that solutions of the reduced equations relate to solutions of the original equations. Moreover, at least for atmospheric flows it is not a completely consensus agreement that there is a natural scale separation. We refer to \cite[Section 1.4]{Klein_2010} where this is briefly discussed. This is quite a natural point of discussion as it is still mathematically a partially unsolved problem, how far one has to track small-scale energy cascades to large scales in geophysical fluid dynamics.
                \end{itemize}

Using characteristic scales can be an intuitive and straightforward approach for a model reduction, provided that one has the knowledge of the existence of such characteristic scales. However, there are oftentimes no direct estimates for approximation errors, so that it can be unclear whether there is a rigorous connection between the original and the reduced models. Invariant manifold reductions, which we are going to discuss next, are among the techniques for which this connection can indeed be described rigorously.

\section{Invariant Manifolds}
\label{sec:invmfld}

Invariant manifolds are a key element in the toolbox of dynamical systems theory. Yet, although different approaches have been used successfully and rigorously in several applications, climate dynamics poses additional intricate problems. This pertains crucially to the level of detail of a reduction one aims for. If one is just interested in the behaviour near a single equilibrium point, stable/unstable and centre manifold theory applies to a wide variety of differential equations. centre manifolds~\cite{Carr} are particularly useful as they isolate the slowest/neutral modes and provide an efficient dimension reduction. At the next level of detail, one can ask for larger invariant manifolds, which contain non-trivial slow dynamics. This brings us to the realm of slow manifolds~\cite{Fenichel_1979,Jones,Kuehn_2015}, which are still generically local and focus on finite time dynamics. Of course, one could even be more ambitious and aim for a global reduction to a single effective low-dimensional dynamical system, where manifolds are usually inertial manifolds~\cite{Robinson1,Temam_1997,KuehnBook1}. Here we present these three viewpoints in reverse order, going from the most global view of inertial manifold, to the intermediate and flexible compromise of slow manifolds, and finally we end describing classical local centre manifold theory near a non-hyperbolic equilibrium.

\subsection{Inertial Manifold Reduction}
\label{sec:Inertial}

The concept of an inertial manifold is classical and aims for a global reduction of a differential equation to a simple, hopefully low-dimensional, set of ordinary differential equations. The most classical setup for inertial manifolds occurs in the context of dissipative evolution equations of the form
\begin{align}\label{Eq:SemilinearEvolutionEquation}
	\partial_t u + Au = f(u),\quad u(0)=u_0.
\end{align}
In this scenario, $A\colon H\supset D(A)\to H$ is usually assumed to be a self-adoint operator with compact inverse on a Hilbert space $H$ and with domain $D(A)$. Moreover, $f\colon H\to H$ is a Lipschitz continuous nonlinearity. Technically, one may allow Lipschitz continuous functions $f\colon H\to D(A^{\beta})$ for $\beta>-1$, see for example \cite[Section 2.5]{Zelik_2014}, and this case is also important if the nonlinearity involves derivatives of $u$. But for the sake of simplicity, we refrain from giving a definition of $D(A^{\beta})$ and work with the simpler case $f\colon H\to H$ instead. The dissipativity of \eqref{Eq:SemilinearEvolutionEquation} is understood in the sense that there are positive constants $C,\alpha>0$ and a nondecreasing function $N\colon[0,\infty)\to [0,\infty)$ such that solutions of \eqref{Eq:SemilinearEvolutionEquation} satisfy the estimate
\[
	\|u(t)\|_{H}\leq N(\|u_0\|_H) e^{-\alpha t} + C\quad (t\geq0).
\]
An inertial manifold $\mathcal{M}$ is a smooth manifold in $H$ such that
\begin{enumerate}[(i)]
	\item\label{InertialManifold:List:FiniteD} $\mathcal{M}$ is finite-dimensional,
	\item\label{InertialManifold:List:Invariance} $\mathcal{M}$ is invariant, i.e., it holds that $T(t)\mathcal{M}\subset\mathcal{M}$, where $(T(t))_{t\geq0}$ denotes the semiflow generated by \eqref{Eq:SemilinearEvolutionEquation}.
	\item\label{InertialManifold:List:ExponentialTracking} $\mathcal{M}$ possesses the exponential tracking property, i.e., there are constants $C,\alpha>0$ such that for all $u_0\in H$ there exists a $v_0\in\mathcal{M}$ such that the corresponding solutions $u$ and $v$ satisfy
		\[
			\|u(t)-v(t)\|_{H}\leq C e^{-\alpha t} \|u_0-v_0\|_{H}\quad(t\geq0).
		\]
\end{enumerate}
Sometimes, especially in earlier works on slow manifold, \eqref{InertialManifold:List:ExponentialTracking} is replaced by the weaker condition
\begin{itemize}
	\item[(iii')] For all $u_0$ there are constants $c_1,c_2>0$ such that the solution $u$ satisfies $\operatorname{dist}(u(t),\mathcal{M})\leq c_1 e^{-c_2 t}$ for all $t\geq0$ (cf. for example \cite{Temam_1990}).
\end{itemize}
However, it turns out that in most relevant scenarios in which inertial manifolds can be shown to exist, they already satisfy \eqref{InertialManifold:List:ExponentialTracking}.

\smallskip

 {\bf{Idea:}} The invariance ensures that an inertial manifold is closed under the dynamics and that the semiflow generated by \eqref{Eq:SemilinearEvolutionEquation} can also be considered just on the inertial manifold instead of the whole underlying space as a phase space. Due to the exponential tracking property, one does not lose any long-term information about the dynamics of \eqref{Eq:SemilinearEvolutionEquation} by this restriction. Nonetheless, since an inertial manifold is finite-dimensional, the restriction of the dynamics to it is expected to yield a much simpler dynamical system than the one corresponding to \eqref{Eq:SemilinearEvolutionEquation}.
 
\smallskip

 {\bf{Outline of the Reduction Procedure:}} Generally, inertial manifolds are known to exist if $A$ satisfies a certain spectral gap condition in relation to $f$, see for example \cite[Theorem 2.1]{Zelik_2014}. More precisely, one usually assumes that the spectrum of $A$ consists of a non-decreasing sequence of positive eigenvalues $(\lambda_k)_{k\in\N}\subset(0,\infty)$. Let $(e_k)_{k\in\N}$ be the corresponding sequence of orthonormal eigenvectors and let $L_f$ be the Lipschitz constant of $f$. If there is an $N\in\N$ such that 
 \begin{align}\label{Eq:SpectralGap}
 	\lambda_{N+1}-\lambda_N>2L_f,
 \end{align}
 then there is an $N$-dimensional inertial manifold $\mathcal{M}$ given as a graph of a Lipschitz-continuous mapping $h\colon H_+\to H_-$, i.e.,
 \[
	 \mathcal{M}=\{u_++h(u_+):u_+\in H_+\}.
 \]
  Here, $H_+:=\operatorname{span}\{e_1,\ldots,e_N\}$ and $H_-$ is the space generated by the orthonormal system $(e_k)_{k\in\N,k>N}$. Let $\operatorname{pr}_{H_+}\colon H\to H_+$ be the canonical projection onto $H_-$. Similar to the procedure for slow manifolds, one may then reduce \eqref{Eq:SemilinearEvolutionEquation} to
  \begin{align*}
  	\partial_tu_++Au_+=\operatorname{pr}_{H_+}f(u_++h(u_+)),\quad u=u_++h(u_+).
  \end{align*}
  This dynamics of this system are now only finite-dimensional.

\smallskip

 {\bf{Example:}} After the concept of an inertial manifold was introduced, there were attempts to prove the existence of inertial manifolds for 2-dimensional Navier-Stokes equations in different geometries, see for example \cite{Kwak_1992,Temam_Wang_1993}. However, there is a gap in the proof in these early attempts so that the existence of inertial manifolds for the Navier-Stokes equations is still an open problem. We refer to \cite{Kostianko_Zelik_2021} in which the techniques from \cite{Kwak_1992,Temam_Wang_1993} are discussed. Since proving the existence of an inertial manifold for Navier-Stokes equations is a difficult problem even in the $2$-dimensional case, inertial manifolds for several regularizations have been considered in recent years, such as Leray $\alpha$-models \cite{AbuHamed_Guo_Titi_2015,Kostianko_2018,Li_Sun_2020}, hyperviscous Navier-Stokes equations \cite{Avrin_2008,Gal_Guo_2018},\cite[Chapter IX.4.3]{Temam_1997} or combinations of both \cite{Kostiankio_Li_Sun_Zelik_2020}. Let us give more details on the hyperviscous case. The underlying equation on the torus is given by
  \begin{align}
  \begin{aligned}\label{Eq:InertialManifold:HNS}
      \partial_t u + \nu(-\Delta)^{l}u +(u\cdot\nabla)u+\nabla p &= g,\quad (t,x)\in \R_+\times \mathbb{T}^n,\\
      \operatorname{\div} u = 0,\quad u\vert_{t=0}&=u_0,
     \end{aligned}
  \end{align}
  where $u$ and $p$ are the unknown functions, $\nu$ is a given parameter, $g$ a given forcing and $u_0$ a given initial value. The parameter $l$ equals to $1$ for the classical Navier-Stokes equation and is taken larger than $1$ in the hyperviscous case. The existence of inertial manifolds however has only been shown for an even smaller set of exponents. But let us first bring \eqref{Eq:InertialManifold:HNS} into the abstract form \eqref{Eq:SemilinearEvolutionEquation}. One may hide the incompressibily condition $\operatorname{\div} u = 0$ as well as the pressure $p$ by projecting the equation to divergence free vector fields using the Helmholtz projection
  \[
    P\colon L_2(\mathbb{T}^n;\R^n)\to L^{\sigma}_2(\mathbb{T}^n;\R^n),
  \]
  where $L^{\sigma}_2(\mathbb{T}^n;\R^n)$ denotes the subspace of $L_2(\mathbb{T}^n;\R^n)$ consisting of the divergence free vector fields. Applying this projection to \eqref{Eq:InertialManifold:HNS} yields
  \begin{align}\label{Eq:InertialManifold:HNS_Projected}
       \partial_t u + \nu (-\Delta)^{l}u +P (u\cdot\nabla)u &= P g,\quad (t,x)\in \R_+\times \mathbb{T}^n.
  \end{align}
  Here, we used that $u$ is divergence free and that $P(-\Delta)^{l} u = (-\Delta)^{l} P u = (-\Delta)^{l}u$. If we now choose 
  \begin{align*}
   & H=L_2^{\sigma}(\mathbb{T}^n;\R^n),\quad D(A):=\{u\in H: (-\Delta)^{\beta}u\in H\},\\
    &A\colon H\supset D(A)\to H,\,u\mapsto (-\Delta)^{l}u,\\
   & f(u)= -P(u\cdot\nabla) u + Pg,
  \end{align*}
  then \eqref{Eq:InertialManifold:HNS_Projected} almost is of the form \eqref{Eq:SemilinearEvolutionEquation}. The only additional step is that one has to carry out a cutoff procedure for $f$ to become Lipschitz continuous. After such a cutoff procedure, one can verify the existence of an inertial manifold for \eqref{Eq:InertialManifold:HNS_Projected} if $l>\frac{3}{2}$ using a spectral gap condition of the form \eqref{Eq_GeneralSpectralGap} which we introduce later. Indeed, the eigenvalues of $(-\Delta)^l$ are of the form $ (|k|_2^{2})^l$ $(k\in\Z^n)$ and therefore, every eigenvalue will be of the form $m^l$ for a natural number $m\in\N_0$. Moreover, each of these eigenvalues has finite multiplicity so that for every $m\in\N_0$ there is an $N\in\N$ such that $\lambda_{N+1}-\lambda_N=(m+1)^l-m^l\geq l m^{l-1}$. For the denominator in  \eqref{Eq_GeneralSpectralGap}, we can choose $\beta=-\frac{1}{l}$, since $(-\Delta)^l$ is of order $2l$ and the nonlinear advection term is of order $1$. Hence, we have $\lambda_{N+1}^{1/2l}+\lambda_N^{1/2l}\leq 2(m+1)^{1/2}$ so that we obtain
  \[
    	\frac{\lambda_{N+1}-\lambda_N}{\lambda_{N+1}^{-\beta/2}+\lambda_{N}^{-\beta/2}}\gtrsim m^{l-\tfrac{3}{2}}\to \infty 
  \]
  as $m\to\infty$ or equivalently $N\to\infty$ if $l>3/2$. 
  This shows that there is an $N\in\N$ depending on the spatial dimension, such that there is an inertial manifold for the hyperviscous Navier-Stokes equation with $l>3/2$ which is given as a graph over the first $N$ eigenfunction of $(-\Delta)^l$ in $L_2^{\sigma}(\mathbb{T}^n;\R^n)$. At least for $n\in\{2,3\}$ the existence of an inertial manifold has also been verified in the critical case $l=3/2$, see \cite{Gal_Guo_2018,Kostiankio_Li_Sun_Zelik_2020}. Therein a technique called spatial averaging was used instead of the spectral gap approach.
  
\smallskip

 {\bf{Comments:}} 	\begin{itemize}
 		\item If $f$ is only Lipschitz continuous as a mapping from $H$ to $D(A^{\beta/2})$ for $\beta\in(-2,0]$, then one has to impose
 		\begin{align}\label{Eq_GeneralSpectralGap}
 			\frac{\lambda_{N+1}-\lambda_N}{\lambda_{N+1}^{-\beta/2}+\lambda_{N}^{-\beta/2}}>L
 		\end{align}
 		instead of \eqref{Eq:SpectralGap} as spectral gap condition, see for example \cite[(1.3)]{Kostianko_Zelik_2021}.
 		\item To the best of our knowledge, the concept of an inertial manifold was first explicitly introduced in \cite{Foias_Sell_Temam_1988}. However, it is mentioned in the introduction of \cite{Foias_Sell_Temam_1988} that related ideas have already been used before, see for example \cite[Chapters 6,9]{Henry_1981,Mane_1977,Mora_1983}.
 		\item   In practice, nonlinearities are oftentimes not globally Lipschitz continuous but only locally. However, the asserted dissipativity of the system allows one to use cut-off techniques to circumvent this issue without changing the global long-term dynamics.
 		\item The main advantage of an inertial manifold reduction is that, if inertial manifolds exist, they yield a global reduction of an infinite-dimensional to a finite-dimensional system without losing any asymptotic information of the dynamics. However, the downside is that their existence can only be shown under very restrictive conditions and for many important examples such as Navier-Stokes equations, their existence is still an open problem. Therefore, the number of scenarios in which an inertial manifold reduction can be applied is rather limited. 
 	\end{itemize}
 	
The main idea behind an inertial manifold reduction is that the dissipativity of a system should in some sense eventually lead to finite-dimensional dynamics. The next reduction method we are going to discuss uses different timescales instead of the dissipativity of the underlying equation. For so-called slow manifold reduction, on usually has a system of a fast and a slow variable and one aims to reduce the dynamics to the ones of the slow variable only. Let us now describe this idea in more detail.

\subsection{Slow Manifold Reduction}
\label{sec:Slow}

Slow manifold reductions are used in systems with multiple time scales. In their standard form such systems are given by
\begin{align}
	\begin{aligned}\label{Eq:FastSlow}
		\epsilon \partial_t u^{\epsilon} &= Au^{\epsilon}+f(u^{\epsilon},v^{\epsilon}),\\
		\partial_t v^{\epsilon} &= Bv^{\epsilon}+g(u^{\epsilon},v^{\epsilon}).
	\end{aligned}
\end{align}
Here, $0<\epsilon \ll1$ is a small parameter indicating the presence of different time scales, $A,B$ are linear operators, $f,g$ are given Lipschitz continuous functions and $u^{\epsilon},v^{\epsilon}$ are the unknown functions. Suppose that $u^{\epsilon}$ takes values in the Banach space $X$ and $v^{\epsilon}$ takes values in den Banach space $Y$. If \eqref{Eq:FastSlow} is an ordinary differential equation, then the framework to reduce its dynamics is provided by classical Fenichel-Tikhonov theory \cite{Fenichel_1979,Tikhonov}, which turned out to be of great importance in multiple time scale dynamics. For partial differential equations however, the theory is much less complete~\cite{Bates_Lu_Zeng_1998}. While slow manifolds have been constructed for certain examples over the years, a first step towards a rigorous justification of slow manifold reductions for abstract partial differential equations of the form~\eqref{Eq:FastSlow} was taken in \cite{Hummel_Kuehn_2020,EngelHummelKuehn}.  

\smallskip

 {\bf{Idea:}} Under appropriate conditions, the fast variable $u^{\epsilon}$ quickly approaches a state in which its dynamics are only determined by the slow variable. More precisely, it gets attracted by a slow manifold $S_{\epsilon}$, a subset of the phase space which is given as a graph over the slow variable space and on which the dynamics of the full system evolve on the slow time scale. The fact that trajectories are attracted by the slow manifold is taken as a justification that one only has to consider the flow on the slow manifold. And since the slow manifold is given as a graph, one can reduce the system to a self-contained, potentially lower-dimensional equation which only depends on the slow variable.
 
\smallskip

 {\bf{Outline of the Reduction Procedure:}} 	The following procedure only works under certain conditions. Below, we will elaborate on them a little bit more, but for the moment we just assume that we can carry out the following steps and that all the objects we use do exist.\\
 One starts with a fast-slow system of the form \eqref{Eq:FastSlow}. One splits the slow variable space $Y=Y^{\zeta}_F\oplus Y^{\zeta}_S$ into a $B$-invariant part $Y^{\zeta}_F$ which decays quickly under the semiflow $(e^{tB})_{t\geq0}$ generated by the linear part $B$ in the slow equation, and a $B$-invariant part $Y^{\zeta}_S$ on which the operators $(e^{tB})_{t\geq0}$ are invertible and their inverses $(e^{-tB})_{t\geq0}$ on $Y^{\zeta}_S$ do not grow quickly. Which rates of decay or growth are considered as 'quickly' is controlled by a parameter $\zeta$. If $\zeta$ gets smaller, then $Y^{\zeta}_F$ gets smaller and $Y^{\zeta}_S$ gets larger. Now, if $\epsilon,\zeta>0$ are small enough and satisfy $\epsilon<c\zeta$ for a certain $c\in(0,1)$, then there is a mapping $h^{\epsilon,\zeta}=(h^{\epsilon,\zeta}_X,h^{\epsilon,\zeta}_{Y_F^{\zeta}})\colon Y_S^{\zeta}\to X\times Y_F^{\zeta}$ with the following properties:
 \begin{enumerate}[(i)]
	\item\label{SlowManifold:List:Invariance} The so-called slow manifold $S_{\epsilon,\zeta}:=\{(h^{\epsilon,\zeta}_X(v),h^{\epsilon,\zeta}_{Y_F^{\zeta}}(v)+v): v\in Y_S^{\zeta}\cap D(B)\} $ is invariant under the semiflow generated by \eqref{Eq:FastSlow}. Here, $D(B)$ denotes the domain of $B$.
	\item\label{SlowManifold:List:Attracting} The slow manifold $S_{\epsilon,\zeta}$ attracts all solutions of \eqref{Eq:FastSlow} exponentially quickly.
	\item\label{SlowManifold:List:DistToCrit} The slow manifold is close to the critical manifold $S_0:=\{(x,y):Ax+f(x,y)=0\}$, where $(x,y)$ is taken such that $Ax+f(x,y)$ is well-defined.
		\item\label{SlowManifold:List:SlowSubsystem} The solutions of $\eqref{Eq:FastSlow}$ on the slow manifold are close to solutions of
		 \begin{align}
	\begin{aligned}\label{Eq:ReducedSlowSubsystem}
		0 &= f(u^{0,\zeta},v^{0,\zeta}),\\
		0 &= \operatorname{pr}_{Y_F^{\zeta}} v^{0,\zeta},\\
		\partial_t v^{0,\zeta} &= Bv^{0,\zeta}+\operatorname{pr}_{Y_S^{\zeta}}g(u^{0,\zeta},v^{0,\zeta}).
	\end{aligned}
\end{align}
Here, $\operatorname{pr}_{Y_F^{\zeta}}\colon Y\to Y_F^{\zeta}$ and $\operatorname{pr}_{Y_S^{\zeta}}\colon Y\to Y_S^{\zeta}$ denote the canonical projections according to the splitting $Y=Y_F^{\zeta}\oplus Y_S^{\zeta}$
 \end{enumerate}
 The properties \eqref{SlowManifold:List:Invariance} and \eqref{SlowManifold:List:Attracting} can be taken as a justification that one may reduce \eqref{Eq:FastSlow} to 
 \begin{align*}
 	u^{\epsilon,\zeta}(t)&=h^{\epsilon,\zeta}_{X}(v_S^{\epsilon,\zeta}(t)),\\
 	v_F^{\epsilon,\zeta}(t)&=h^{\epsilon,\zeta}_{Y_F^{\zeta}}(v_S^{\epsilon,\zeta}(t)),\\
 	\partial_t v_S^{\epsilon,\zeta}(t) &= Bv_S^{\epsilon,\zeta}(t)+\operatorname{pr}_{Y_S^{\zeta}}g(h^{\epsilon,\zeta}_{X}(v_S^{\epsilon,\zeta}(t)),h^{\epsilon,\zeta}_{Y_F^{\zeta}}(v_S^{\epsilon,\zeta}(t))+v_S^{\epsilon,\zeta}(t)).
 \end{align*}
 The dynamical part of this system is only determined by the last equation and depends only on the slow part of the slow variable. It may even be finite-dimensional in many cases, even if both $X$ and $Y$ are infinite-dimensional. The other components are then uniquely determined by $v_S^{\epsilon,\zeta}(t)$. Moreover, the properties \eqref{SlowManifold:List:Attracting} and \eqref{SlowManifold:List:SlowSubsystem} give us information about the location of the slow manifold and the flow on it, even if it is difficult to construct it explicitely in certain scenarios. \\
 The proximity of the slow manifold to the critical manifold also motivates a simpler reduction of \eqref{Eq:FastSlow}, namely the one to the slow subsystem given by
 \begin{align}
	\begin{aligned}\label{Eq:SlowSubsystem}
		0 &= f(u^{0},v^{0}),\\
		\partial_t v^{0} &= g(u^{0},v^{0}).
	\end{aligned}
\end{align}
This reduction is less accurate in the sense that the flow generated by \eqref{Eq:SlowSubsystem} is not contained in the original one generated by \eqref{Eq:FastSlow}. However, the slow subsystem is still a reasonable approximation of \eqref{Eq:FastSlow}. In addition, this reduction has more advantages. Most importantly, it does not need the existence of the above splitting $Y=Y^{\zeta}_F\oplus Y^{\zeta}_S$ of the slow variable space, which is a rather restrictive assumption. The reduction to the slow subsystem is therefore applicable in many more situations, even if slow manifolds do not exist.

\smallskip

 {\bf{Example:}}  A well-known example of a slow manifold reduction is the one for the different variants of the Stommel model. Here, we take the Stommel-Cessi model given by the fast-slow system
\begin{align}\label{Eq:SlowManifold:Stommel}\begin{aligned}
    \epsilon\dot{u} &= -(u-1) - \epsilon u [1+7.5(u-v)^2],\\
    \dot{v} &= \mu - v[1+7.5(u-v)^2],
\end{aligned}\end{align}
see \cite{Cessi_1994,Stommel_1961} and \cite[(6.2.27)]{Berglund_Gentz_2006}. Here, $u$ denotes the temperature difference and $v$ the salinity difference between a higher and a lower latitude box in the North Atlantic. The parameter $\mu$ is proportional to the freshwater flux and can mathematically be seen as a bifurcation parameter for the dynamics of the thermoaline circulation in the North Atlantic. By Fenichel theory, there is a slow manifold which is a perturbation of order $\epsilon$ of the critical manifold given by the constant function $h^0(v)\equiv1$. One may thus approximate \eqref{Eq:SlowManifold:Stommel} by
\begin{align}\label{Eq:SlowManifold:Stommel_Reduced}
    \dot{v}=\mu-v[1+7.5(1-v)^2].
\end{align}
After this reduction, the bifurcation analysis of this equation boils down to studying the roots of the third order polynomial on the right-hand side. One can verify that it undergoes two saddle-node bifurcations at
\[
    \mu_{\pm}=\frac{1}{9}\left(11\pm3\cdot\sqrt{\frac{3}{5}}\right).
\]
If $\mu<\mu_-$ or $\mu>\mu_+$ there is only one stable equilibrium which corresponds to a low salinity difference or high salinity difference, respectively. If $\mu\in(\mu_-,\mu_+)$, then there are two stable states the system may be in.\\
A similar reduction can be carried out in certain situations if we add a spatial component. For example, the system
\begin{align}\label{Eq:SlowManifold:Diffusive_Stommel}\begin{aligned}
    \epsilon\partial_{t} u(t,x) &= \partial_{xx}u(t,x) -(u(t,x)-1) - \epsilon u(t,x) [1+7.5(u(t,x)-v(t,x))^2],\\
    \partial_{t} v(t,x) &= \partial_{xx}v(t,x) \mu - v(t,x)[1+7.5(u(t,x)-v(t,x))^2],
\end{aligned}\end{align}
on the $1$-dimensional torus, i.e., with $x\in\mathbb{T}$ and periodic boundary conditions, can be approximated by the equation
\begin{align}\label{Eq:SlowManifold:Diffusive_Stommel_Reduced}
     \partial_{t} v(t,x) &= \partial_{xx}v(t,x) \mu - v(t,x)[1+7.5(1-v(t,x))^2],\quad u(t,x)\equiv1.
\end{align}
One may even reduce this equation further to a system of ordinary differential equations by truncating to a finite set of Fourier modes. We refer to \cite[Section 6.1]{Hummel_Kuehn_2020} for more details on the reduction.

\smallskip

 {\bf{Comments:}} \begin{itemize}
                 \item Let us now elaborate a bit more on the necessary assumptions: If \eqref{Eq:FastSlow} is an ordinary differential equation, then it is comparably easy to put into an abstract framework. The functions $f$ and $g$ are naturally given by the application and act on $X\times Y= \R^m\times\R^n$ as a phase space, i.e., $u^{\epsilon}$ takes values in $X=\R^m$ and $v^{\epsilon}$ takes values in $Y=\R^n$ for some $n,m\in\N$. If the nonlinearities are smooth enough, then the only condition one needs for the reduction procedure is that the linearization of $u\mapsto Au+f(u,v)$ does not have eigenvalues on the imaginary axis for every $v\in \R^n$. The splitting of the slow variable space is trivial for ordinary differential equations: One can just take $Y_F^{\zeta}=\{0\}$ and $Y_S^{\zeta}=\R^n$. The slow manifolds $S_{\epsilon,\zeta}$ will then be independent of $\zeta$ so that they only depend on $\epsilon$ as in classical Fenichel theory. 
 \item	However, if \eqref{Eq:FastSlow} is a partial differential equation, it requires more effort to put it into an abstract form which allows one to treat it the setting of \cite{Hummel_Kuehn_2020}. The main reason is that in infinite dimensions, one has many more possible choices of phase spaces and corresponding topologies. The assumptions for the theory in \cite{Hummel_Kuehn_2020} might not be satisfied in some natural looking choice of phase spaces while they may be satisfied in other phase spaces. Moreover, the functions $f$ and $g$ will usually not act on the whole phase space but only on a certain subset. If for example $X=C^{\alpha}(\R)$ is the space of $\alpha$-H\"older continuous functions on $\R$ with some $\alpha\in(0,1)$ and $f(u^{\epsilon},v^{\epsilon})=\Delta u^{\epsilon}$, then one formally needs that $u^{\epsilon}$ is twice differentiable and that its derivatives are $\alpha$-H\"older continuous. One is therefore forced to work with different spaces and different topologies.
 Aside from the technicalities concerning different topologies and spaces, there are two important conditions one has to verify. Firstly, the semigroup generated by $A$ should have a negative growth bound and the Lipschitz constant of $f$ should be small enough so that solution of the fast equation would decay without the forcing from the slow variable. Secondly, one has to find a splitting $Y=Y_F^{\zeta}\oplus Y_S^{\zeta}$ as mentioned above. The space $Y_F^{\zeta}$ contains the parts which decay quickly under the $C_0$--semigroup $(e^{tB})_{t\geq0}$ generated by $B$. On $Y_S^{\zeta}$ the operators $(e^{tB})_{t\geq0}$ are invertible with inverses $(e^{-tB})_{t\geq0}$ these inverses do not grow too quickly as $t\to\infty$. Moreover, the rates of decay and growth of $(e^{tB})_{t\geq0}$ and $(e^{-tB})_{t\geq0}$, respectively, have to be seperated well enough so that it persists even if the nonlinearity $g$ is taken into account. The seperation should even become stronger as $\zeta$ gets smaller for convergence results to hold true. Oftentimes this means that $B$ has to have spectral gaps which get arbitrarily large as one moves towards the left in the complex plane. Unfortunately, this assumption is very restrictive. If for example $B=\Delta$ is the Laplacian on the $n$-dimensional torus $\mathbb{T}^n$, then the distance between the eigenvalues does not get arbitrarily large unless $n\in\{1,2\}$. In the case of only small spectral gaps or even the absence of spectral gaps, one is still forced to rely on the less accurate reduction to the slow subsystem.
 \item Nonlinearities in many applications are actually not Lipschitz continuous, but only locally Lipschitz-continuous. Nonetheless, they can be made Lipschitz continuous by a standard cut-off procedure. In such cases however, the slow manifolds also only have a local meaning instead of a global one.
 \item A slow manifold reduction allows one to remove the fast variable $u^{\epsilon}$ from the dynamics. Moreover, this method does not only reduce the dimension of the dynamics, but also provides a nice geometric interpretation of the reduction, since the reduced system is contained in an attracting invariant manifold of the original system. However, the assumptions of a slow manifold reduction are not directly satisfied at bifurcation points of the fast variable. Thus, many interesting pattern forming problems are not accessible through a slow manifold reduction alone and instead, one has to use additional desingularization techniques such as the blow-up method, see for example \cite{Engel_Kuehn_2020,Krupa_Szmolyan_2001a,Krupa_Szmolyan_2001b}.
 \item Difficulties also arise if the equation for the slow variable is given by a partial differential equations. In this case, the conditions under which a slow manifold reduction can be rigorously justifies are as restrictive as the ones for an inertial manifold reduction. However, one may still approximate trajectories of the fast-slow system by the ones of the slow subsystem under less restrictive assumptions also in the infinite-dimensional case, see \cite{Hummel_Kuehn_2020}.
                \end{itemize}

While inertial and slow manifold reductions aim to describe the dynamics of a system in a larger part of the phase space, we are now going to discuss a reduction method which only aims to reduce the dynamics around a non-hyperbolic equilibrium: the so-called centre manifold reduction.

\subsection{Centre Manifold Reduction}
\label{Sec:centre}

Another important classical reduction method using invariant manifolds is the centre manifold reduction. It is used to simplify the dynamics of a system of the form
\begin{align}\label{Eq:centreManifold}
	\dot{x}=f(x)
\end{align}
 locally around a non-hyperbolic fixed point $x_0\in X$ with $X$ being a Banach space. A common extension is to consider parametrized bifurcation problems of the form
 \begin{align}\label{Eq:centreManifoldParams}
	\dot{x}=f(x,\lambda)
\end{align}
with $\lambda$ a finite dimensional parameter set \cite{kuznetsov2013elements}.
 \smallskip

 {\bf{Idea:}}  A centre manifold consists of those trajectories around a fixed point, which are neither attracted nor repelled by it at an exponential rate. Roughly speaking, the long-term behaviour of the former trajectories shrinks down to the equilibrium point itself, while the latter ones are pushed away from the fixed point so that they are not important for the local dynamics. The trajectories on the centre manifold in turn persist and stay close to the fixed point. It is therefore sufficient to restrict the system to the centre manifold in order to describe its effective dynamics around the fixed point.
  
\smallskip

 {\bf{Outline of the Reduction Procedure:}} Assume that $f$ is smooth enough. As already mentioned above, centre manifolds are constructed around a non-hyperbolic steady state, i.e., around an $x_0\in X$ such that $f(x_0)=0$ and that the spectrum of the linearization $Df(x_0)$ splits into three parts
 \[
	\sigma( Df(x_0) ) = \sigma^s\,\,\cup\sigma^c\,\cup\,\sigma^s
 \]
 with $\sigma_c\neq\emptyset$, where 
 \begin{align*}
	\sigma^s&:=\{\lambda\in \sigma( Df(x_0) ) : \operatorname{Re}\lambda<0\} ,\\
	\sigma^c&:=\{\lambda\in \sigma( Df(x_0) ) : \operatorname{Re}\lambda=0\} ,\quad\text{and}\\
	\sigma^u&:=\{\lambda\in \sigma( Df(x_0) ) : \operatorname{Re}\lambda>0\}
 \end{align*}
 denote the stable, neutral (centre), and unstable part of the spectrum, respectively. In infinite dimensions one has to impose additional conditions which partly depend on the desired properties of the centre manifold. One possibility, which is carried out in \cite{Bates_Jones_1989}, is to assume that:
 \begin{enumerate}[(i)]
 	\item The linearization $Df(x_0)$ generates a $C_0$-semigroup $(T(t))_{t\geq0}\subset\mathcal{B}(X)$.
 	\item The three parts of the spectrum $\sigma^s,\sigma^c$ and $\sigma^u$ are closed and open, and there are corresponding spectral projections $\operatorname{pr}_s\colon X\to X^s$, $\operatorname{pr}_c\colon X\to X^c$ and $\operatorname{pr}_u\colon X\to X^s$, respectively.
 	\item $\operatorname{dim} X^c<\infty$, $\operatorname{dim} X^u<\infty$.
 	\item There are constants $M,c>0$ such that $\|T(t)\vert_{X^s}\|_{\mathcal{B}(X^s)}\leq Me^{-ct}$ holds for all $t\geq0$.
 \end{enumerate}
 Other assumptions are for example discussed in \cite{Vanderbauwhede_Iooss_1992}. Note that one can adjust the assumptions so that it is not necessary to assume that the centre space $X^c$ is finite-dimensional, as the conservative case in \cite{Bates_Jones_1989} shows. Even though the assumptions might differ slightly depending on the particular situation, one locally constructs a centre manifold $\mathcal{C}$ as the graph of a mapping $h\colon \{x_0\}+U\to X^s\oplus X^u$ over the centre subspace, i.e.,
 \[
 	\mathcal{C}=\{h(x_0+x)+x_0+x^c:x\in U\},
 \]
 where $U\subset X^c$ is a certain neighborhood of $0$ in $X^c$. This is usually done by one of the standard approaches such as a Lyapunov-Perron or a graph transform argument. Similar to the slow and inertial manifold reductions, one can now reduce \eqref{Eq:centreManifold} to a dynamical system in $U\subset X^c$, more precisely to
 \[
 	\dot{x}^c=\operatorname{pr}_cf(x^c+h(x^c)), \quad x=x^c+h(x^c)
 \]
 with initial conditions in $\{x_0\}+U$.
   
\smallskip

 {\bf{Example:}} A slow manifold reduction can oftentimes also be seen as a specific case of a centre manifold reduction. To this end, we rescale time by $\tau=\epsilon t$ so that \eqref{Eq:FastSlow} turns into
 			\begin{align}\label{Eq:FastSlow_FastTime}
 			    \begin{aligned}
 				\partial_{\tau} u^{\epsilon} &= Au^{\epsilon} + f(u^{\epsilon},v^{\epsilon}),\\
 				\partial_{\tau} v^{\epsilon} &= \epsilon Bv^{\epsilon} + \epsilon g(u^{\epsilon}, v^{\epsilon}).
 				\end{aligned}
 			\end{align}
 			Now one adds the new equation
 			\[
 				\partial_{\tau} w^{\epsilon} = 0,\quad w^{\epsilon}(0)=\sqrt{\epsilon}
 			\]
 			to the system. This yields the new system
 			 			\begin{align}\label{Eq:FastSlow_with_dummy}\begin{aligned}
 				\partial_{\tau} u^{\epsilon} &= Au^{\epsilon} + f(u^{\epsilon},v^{\epsilon}),\\
 				\partial_{\tau} v^{\epsilon} &= (w^{\epsilon})^2 Bv^{\epsilon} + (w^{\epsilon})^2 g(u^{\epsilon}, v^{\epsilon}),\\
 					\partial_{\tau} w^{\epsilon} &= 0,\quad w^{\epsilon}(0)=\sqrt{\epsilon}
 					\end{aligned}
 			\end{align}
 			Setting $\epsilon=0$ leads to the fast subsystem
 			\begin{align*}
 			\partial_{\tau} u^{\epsilon} &= Au^{\epsilon} + f(u^{\epsilon},v^{\epsilon}),\\
 			\partial_{\tau} v^{\epsilon} &= 0,\\
 			\partial_{\tau} w^{\epsilon} &= 0.
 			\end{align*}
 			Let now $(u,v)\in S_0$ be a point on the critical manifold $S_0$ of \eqref{Eq:FastSlow_FastTime}. If the linearization of $A+f$ has no spectrum on the imaginary axis, then $(u,v,0)$ is a non-hyperbolic equilibrium point of \eqref{Eq:FastSlow_with_dummy} with centre subspace $\{0\}\times Y\times \R$. Hence, if $\dim Y<\infty$ one can use centre manifold theory to obtain a neighbourhood $U\subset Y\times\R$ of $(v,0)$ and a centre manifold $S$ of \eqref{Eq:FastSlow_FastTime} containing $(u,v,0)$ which is given as a graph over $U$. For a fixed $\epsilon>0$, the set $S_{\epsilon}:=S\cap\{(u,v,w)\in X\times Y\times \R: w=\epsilon\}$ now defines a slow manifold. Therefore, one could again take the Stommel model \eqref{Eq:SlowManifold:Stommel} as an example for a centre manifold reduction, where the slow manifold is the perturbation of the critical manifold that consists of steady states for the fast subsystem. 
  
\smallskip

 {\bf{Comments:}} \begin{itemize}
 			\item One of the advantages of centre manifold theory is that is it well-developed and widely used, also in the infinite-dimensional or stochastic setting. We refer to \cite{Carr_1981,Escher_Simonett_1998,Haragus_Iooss_2011,Kuehn_Neamtu_2020,Simonett_1995,Vanderbauwhede_Iooss_1992}.
 			\item Compared to an inertial manifold reduction, a centre manifold reduction usually has the drawback that it is not a global but only a local reduction method. Hence, the reduction is only valid in a certain part of the phase space.
 		\end{itemize}

The main benefits of invariant manifold reductions are that they can be rigorously described with explicit estimates for the approximation error, and that they give a nice geometric interpretation of the reduction. However, the latter is not always possible, especially in systems with strong oscillations. In such cases, one can try to use other classical reduction techniques called averaging and homogenization, which we are going to discuss in the next section.

\section{Averaging and Homogenization}
\label{Sec:AvHom}

Averaging and Homogenization are important and similar methods of reducing systems with strong oscillations to averaged/homogenized systems without these oscillations. Both methods are particularly useful for the derivation of stochastic climate models, as they can in certain situations characterize in which sense solutions of \eqref{HasselmannForm} are approximated by solutions of an equation with random forcing. In other words, the intuitive idea that small-scale nonlinear self-interactions can effectively be replaced by a random forcing can be made rigorous. There are many works treating the mathematical ideas behind this theory, see for example \cite{Duan_Wang_2014,Gottwald_Melbourne_2013,Khasminskii_1963,Melbourne_Stuart_2011,Pavliotis_Stuart_2008}. Other works such as \cite{Arnold_2001,Franzke_et_al_2015,Kifer_2001,Roedenbeck_Beck_Kantz_2001} explain them with their importance for climate models in mind. Explicit applications to climate models are mainly given by Majda and co-authors, see for example \cite{Majda_Franzke_Crommelin_2009,Majda_Franzke_Khouider_2008,Majda_Timofeyev_VandenEijnden_2001} and references therein.

We are not aware of a coherent abstract framework that includes both, or a classification of the problems which can be treated by the methods as well as a general reduction procedure which applies to all relevant examples. In fact, the literature is not even consistent about the classification of concrete examples: the reduction of a parabolic partial differential equation with small diffusion and strongly varying coefficients for the advection terms such as the one in \cite[Section 14, Section 21]{Pavliotis_Stuart_2008} is referred to as averaging by some authors and as homogenization by others. Sometimes, both names seem even to be used interchangeably, see for example \cite{Campillo_Kleptsyna_Piatnitski_2001,Diop_Pardoux_2004}. The reason is that there are different conceptions of what averaging and homogenization should mean.

One possible definition is given in \cite[Section 1.3]{Pavliotis_Stuart_2008}. Therein, the authors write that averaging and homogenization can be applied to perturbations of equations of the form
\[
	L^{\epsilon}u^{\epsilon}=f
\]
or
\[
	\partial_t u^{\epsilon} = L^{\epsilon}u^{\epsilon}
\]
where the linear operator $L^{\epsilon}$ has the form
\begin{align}\label{Eq:FirstOrderPerturbation}
	L^{\epsilon}=\epsilon^{-1}L_0+L_1
\end{align}
or
\begin{align}\label{Eq:SecondOrderPerturbation}
	L^{\epsilon}=\epsilon^{-2}L_0+\epsilon^{-1}L_1+L_2.
\end{align}
Usually, it is assumed that $L_0$ has a non-trivial null space and it is argued that the interesting behaviour takes place on this null space. The techniques treating problems with \eqref{Eq:FirstOrderPerturbation} are referred to as averaging, while those treating problems with \eqref{Eq:SecondOrderPerturbation} are called homogenization.

Another possibility to distinguish averaging and homogenization is to distinguish between temporal and spatial oscillations. In this case, averaging is usually referred to techniques which simplify temporal effects, while homogenization is used to simplify spatial effects. This distinction was for example made in \cite{Duan_Wang_2014}, even though it was not explicitly stated as a definition. One should note that this is not consistent with the definition above from \cite[Section 1.3]{Pavliotis_Stuart_2008}. For this work, it is not our aim to decide what the best characterisations of averaging and homogenization should be. Instead, we just formulate the idea for both methods at once and afterwards, we describe the two reduction methods in scenarios in which the classification is widely accepted.
  
\smallskip

 {\bf{Idea:}} In systems with different spatial or temporal scales it is not essential to study the details of the behaviour of the system on finer scales. Instead, only the average or homogeneous part of the high-resolution/fast effects is important for the macroscopic behaviour. In order to make this precise, one has to derive an averaged or homogenized version of the system and to show that solutions of the original and the averaged or homogenized system are close to each other in a suitable sense.

\subsection{Averaging}
\label{Sec:Averaging}

One of the standard situations in which averaging is applied is similar to the one in which a slow manifold reduction is used, but with different conditions. Consider the fast-slow system of stochastic differential equations
\begin{align}
	\begin{aligned}\label{Eq:Averaging:FastSlow}
		\epsilon du^{\epsilon} &= f(u^{\epsilon},v^{\epsilon})\,dt+ \sqrt{\epsilon}\sigma_1(u^{\epsilon},v^{\epsilon})\,dW_1(t),\\
		 dv^{\epsilon} &= g(u^{\epsilon},v^{\epsilon})\,dt+\sigma_2(u^{\epsilon},v^{\epsilon})\,dW_2(t),\\
		 u^{\epsilon}(0)&=u_0,\quad v^{\epsilon}(0)=v_0.
	\end{aligned}
\end{align}
Here, $W_1$ is a $k$-dimensional and $W_2$ an $l$-dimensional Brownian motion which are independent and the nonlinearities $f\colon\R^m\times\R^n\to\R^m$, $g\colon \R^m\times\R^n\to\R^n$, $\sigma_1\colon\R^m\times\R^n\to\R^{m\times k}$ and $\sigma_2\colon\R^{m}\times\R^n\to \R^{n\times l}$ are assumed to be Lipschitz continuous. The aim is to reduce \eqref{Eq:Averaging:FastSlow} to an equation of the form
   \begin{align}\label{Eq:Averaging:FastSlow:Averaged}
 dv(t) = \bar{g}(v(t))\,dt+\bar{\sigma}(v(t))\,dW_1(t),\quad v(0)=v_0,
 \end{align}
 where $\bar{g}$ and $\bar{\sigma}$ are certain averaged versions of $g$ and $\sigma$. The solution $v$ of \eqref{Eq:Averaging:FastSlow:Averaged} is supposed to approximate $v^{\epsilon}$. Formally, the equations look similar to the ones which are suitable for a slow manifold reduction. The main difference concerning the assumptions compared to the setting for the slow manifold reduction is not only the stochastic forcing we allow. Also in the deterministic setting the assumptions are less restrictive: While slow manifold reductions are carried out for systems in which the fast subsystem is normally hyperbolic in a suitable sense, averaging treats the case in which the fast variable is assumed to be sufficiently mixing.
  
\smallskip

 {\bf{Outline of the Reduction Procedure:}}  One starts with the fast-slow system \eqref{Eq:Averaging:FastSlow}. There are many slight variants for the main assumption and the precise approximation results, see for example \cite[Chapter II.3]{Skorokhod_1989}, \cite[Chapters 11,12, 17 and 18]{Pavliotis_Stuart_2008} or \cite[Chapter 5]{Duan_Wang_2014}. Here, we briefly summarize the approach from \cite[Chapter II.3]{Skorokhod_1989}. Apart from smoothness and growth assumptions on the nonlinearities $f,g,\sigma_1$ and $\sigma_2$ as well as boundedness assumptions on the fast variable $u^{\epsilon}$ (we refer to \cite[Chapter II.3, Theorem 12]{Skorokhod_1989} for the precise conditions), it is crucial that \eqref{Eq:Averaging:FastSlow} satisfies the following ergodicity condition: For each fixed $v_0\in\R^n$ we consider the equation
 \[
 du(t) = f(u(t),v_0)\,dt+\sigma_1(u(t),v_0)\,dW_1(t),\quad u(0)=u_0.
 \]
  We assume that for all $v_0\in\R^n$ this equation has a solution $u_{u_0,v_0}$ and an invariant ergodic distribution $\mu_{v_0}$ such that for all $r>0$ and all $\phi\in C(\R^m)$ the uniform ergodicity condition
  \[
	\lim_{T\to \infty}\sup_{|u_0|,|v_0|\leq r}\mathbb{E}\bigg|\frac{1}{T}\int_0^T\phi(u_{u_0,v_0}(t))\,dt - \int \phi(z)\,d\mu_{v_0}(z)\bigg|=0
	\]
	is satisfied. In this case, one may define
 \[
	\bar{g}(v):= \int_{\R^n} g(u,v)\,d\mu_{v}(u)
 \]
 and $\bar{\sigma}(v)\in\R^{n\times n}$ such that
 \[
 	\bar{\sigma}(v)^2=\int_{\R^m}\sigma_2(u,v)\sigma_2(u,v)^\top \,d\mu_v(u).
 \]
  Moreover, the process $v^{\epsilon}$ from \eqref{Eq:Averaging:FastSlow} converges in distribution to the process $v$ which is defined as the solution of 
   \[
 dv(t) = \bar{g}(v(t))\,dt+\bar{\sigma}(v(t))\,dW_1(t),\quad v(0)=v_0.
 \]
   
\smallskip

 {\bf{Example:}} Similar techniques have already been applied to climate models. Let us continue with the example from Section~\ref{Sec:Galerkin} that we took from \cite{Majda_Timofeyev_VandenEijnden_2001} where the equations of a barotropic flow on a beta plane with topography and mean flow given by 
\begin{align}
    \begin{aligned}\label{Eq:Averaging:BarotropicFlow1}
    0&=\frac{\partial q}{\partial t}+\nabla^{\bot}\psi\cdot\nabla q + U\frac{\partial q}{\partial x}+\beta\frac{\partial\psi}{\partial x},\\
    q&=\Delta\psi+h,\\
    \frac{dU}{dt}&=\dashint h\frac{\partial\psi}{\partial x}.
    \end{aligned}
 \end{align}
 have been reduced in several steps. Here, $q(x,y,t)$ denotes the small-scale potential vorticity, $U(t)$ the mean flow, $\psi(x,y,t)$ the small-scale stream function, $h(x,y)$ models the topography and $\beta$ the variation of the Coriolis parameter. A couple of reduction procedures that resemble the averaging method presented here are carried out in \cite{Majda_Timofeyev_VandenEijnden_2001}. Let us briefly look at one of the examples, namely the one considered in \cite[Section 7.4]{Majda_Timofeyev_VandenEijnden_2001}. After a truncation to finitely many Fourier modes as in Section~\ref{Sec:Galerkin} the mean flow $U$ is declared as the only climate variable and all other variables are treated as unresolved variables. The nonlinear self-interactions of the unresolved variables are replaced by a stochastic forcing with drift and since one is looking for real-valued solutions, one may impose $w_k=w_{-k}^*$, $(k\in\Z^2)$, for any of the Fourier coefficients coming from \eqref{Eq:Averaging:BarotropicFlow1}. This way, the equations
                 \begin{align}
                    \begin{aligned}\label{Eq:Averaging:BarotropicFlow2}
                     \txtd U &= \frac{2}{\epsilon}  \sum_{k=(k_x,k_y)\in\bar{\sigma}} \operatorname{Im}\left(\frac{k_x\hat{h}_k^*}{|k|}w_k\right)\,\txtd t,\\
                     \txtd w_k& = \frac{i k_x \hat{h}_k}{\epsilon|k|}U\,\txtd t - \frac{i}{\epsilon}\big(k_x U-\frac{k_x\beta}{|k|^2}\big)e_k\,\txtd t\\
                     &\qquad-\frac{1}{\epsilon^2}\gamma_k(w_k-\bar{w}_k)\txtd t+\frac{\sigma_k}{\epsilon}\,\txtd W_k(t)\quad(k=(k_x,k_y)\in \bar{\sigma})
                    \end{aligned}
                 \end{align}
                 are derived. Here, $\hat{h}_k^*$ denotes the complex conjugate of the $k$-th Fourier coefficients of $h$. $\gamma_k,\sigma_k,\bar{w}_k,\epsilon$ are real parameters coming from the replacement of the nonlinear self-interactions by the stochastic forcing. The parameter $\gamma_k$ measures the strength of the drift to the mean $\bar{w}_k$, $\sigma_k$ the strength of the noise and $\epsilon$ denotes the scale separation parameter. Moreover, the $W_k(t)$ are independent Brownian motions and $\bar{\sigma}$ is a subset of $\{k\in\Z^2:|k|^2\leq N\}$ $(N\in\N)$ that contains exactly one of the terms $k$ and $-k$ $(k\in\sigma)$. Except for the $\mathcal{O}(\epsilon^{-1})$-terms in the equations for the unresolved variables $w_k$, \eqref{Eq:Averaging:BarotropicFlow2} already has the form \eqref{Eq:Averaging:FastSlow_Generalized}. And indeed, as was shown in \cite[Theorem 7.9]{Majda_Timofeyev_VandenEijnden_2001} this system can be reduced to
                 \[
                    \txtd U = -\gamma_u(U-\bar{U})\,\txtd t +\sigma_u\,\txtd W(t),
                 \]
                 where
                 \begin{align*}
                   \gamma_u&=2\sum_{k\in\bar{\sigma}}\frac{k_x^2|\hat{h}_k|^2(\hat{h}_k-|k|^2\bar{w}_k)}{|k|^2\gamma_k^2\hat{h}_k},\quad\bar{U}=-\frac{2\beta}{\gamma_u}\sum_{k\in\bar{\sigma}}\frac{k_x^2\bar{w}_k|\hat{h}_k|^2}{|k|^2\gamma_k\hat{h}_k},\\&\qquad\qquad\sigma_u=2\left(\sum_{k\in\bar{\sigma}}\frac{k_x^2\sigma_k^2|\hat{h}_k|^2}{\gamma_k^2}\right)^{1/2}.
                 \end{align*}
    
\smallskip

 {\bf{Comments:}} \begin{itemize}
                 \item Averaging techniques have been further developed in the past few decades and can also be applied to certain classes of stochastic partial differential equations. We refer the reader to \cite[Section 5.2]{Duan_Wang_2014}, where stochastic partial differential equations of the form
                 \begin{align*}
                  \partial_tu^{\epsilon}&=\tfrac{1}{\epsilon}\Delta u^{\epsilon}+\tfrac{1}{\epsilon}f(u^{\epsilon},v^{\epsilon})+\frac{\sigma_1}{\sqrt{\epsilon}}\partial_tW_1,\\
                  \partial_tv^{\epsilon}&=\Delta v^{\epsilon}+g(u^{\epsilon},v^{\epsilon})+\sigma_2\partial_tW_2
                 \end{align*}
               are studied. Here $f,g$ satsify certain Lipschitz conditions and $W_1,W_2$ are two independent trace class Brownian motions in $L_2(\mathcal{O})$. Similar problems have also been treated in \cite{Cerrai_2011,Roberts_Wang_2012}.
               \item Equation \eqref{Eq:Averaging:FastSlow} can also be seen as a special case of 
               \begin{align}
               	\begin{aligned}\label{Eq:Averaging:FastSlow_Generalized}
               		\frac{du^{\epsilon}}{dt}&=\frac{1}{\epsilon^2}f(u^{\epsilon},v^{\epsilon})+\frac{1}{\epsilon}\sigma_1(u^{\epsilon},v^{\epsilon}) \frac{dW_1(t)}{dt},\\
               		\frac{dv^{\epsilon}}{dt}&=\frac{1}{\epsilon}g_1(u^{\epsilon},v^{\epsilon})+g_2(u^{\epsilon},v^{\epsilon})+\sigma_2(u^{\epsilon},v^{\epsilon})\frac{dW_2(t)}{dt},
               	\end{aligned}
               \end{align}
               with $g_1=0$. In \cite[Chapter 11]{Pavliotis_Stuart_2008} an outline is given on how to reduce such a system if the $\sigma_1\sigma_1^\top $ is uniformly positive definite, if  $\mu_v$ has a density and if the centering condition $\int g_1(u,v)\,d\mu_v(u)=0$ holds for all $v$. In such a case, \eqref{Eq:Averaging:FastSlow_Generalized} can also be reduced to an equation of the form
                  \[
 dv(t) = \bar{g}(v(t))\,dt+\bar{\sigma}(v(t))\,dW_1(t),
 \]
 but this time with
 \[
 	\bar{g}(v) = \int g_2(u,v)+(\nabla_v\Phi(u,v))g_1(u,v) d\mu_v(u)
 \]
 and $\bar{\sigma}(v)$ given by the relation 
 \[
 \begin{aligned}
 	\bar{\sigma}(v)\bar{\sigma}(v)^T =& \int \left[\sigma_2(u,v) \sigma_2(u,v)^T+\right.\\
 	& \left. g_1(u,v)\otimes \Phi(u,v)+[g_1(u,v)\otimes \Phi(u,v)]^T \right] d\mu_v(u).
\end{aligned}
 \]
 Here, $\Phi$ is defined as the solution of 
 \[
 	-\langle f,\nabla\rangle \Phi - \frac{1}{2}\operatorname{tr}(\sigma_1^T\sigma_1 D^2_y) \Phi=f_1,\quad \int \Phi(u,v)\,d\mu_v(u) = 0
 \]
 with $D_y^2$ denoting the Hessian matrix and $\operatorname{tr}$ denoting the trace of a matrix. Moreover, for two vectors $x,y\in\R^l$ the tensor product $x\otimes y$ is defined to be the matrix $x\otimes y\in\R^{l\times l}$ such that $(x\otimes y)z=\langle y,z\rangle x$ for all $z\in\R^l$. In \cite[Chapter 18]{Pavliotis_Stuart_2008} the reduction is studied more rigorously, but in a more specific setting. Note however, that the reduction of \eqref{Eq:Averaging:FastSlow_Generalized} is called homogenization instead of averaging in \cite{Pavliotis_Stuart_2008}.
 \item Related methods, which are often referred to as fast wave averaging, have also been applied the geophysical flows. These methods are based on techniques that were developed in \cite{klainerman1981singular,schochet1994fast} for singular limits of quasilinear hyperbolic systems. For example, the systems that are treated in \cite{schochet1994fast} are of the form
 \begin{align}\label{Eq:SymHypSystem}
 \begin{aligned}
    &A_0(\epsilon U^{\epsilon}) \partial_t U^{\epsilon} +\sum_{i=1}^N A_{i}(U^{\epsilon},\epsilon) \partial_{x_i} U^{\epsilon} +\\ 
    & \qquad \sum_{j=1}^M \left(\frac{1}{\epsilon}K_j+D_j(U^{\epsilon},\epsilon)\right) \partial_{y_j}U^{\epsilon} = F(U^{\epsilon},\epsilon),\\
    &\qquad U^{\epsilon}(0,x,y) = U_0(x,y)
    \end{aligned}
 \end{align}
 where $N,M\in\N$, $A_0$ takes values in the set of positive definite matrices and $A_i,K_j,D_j$ take values in the set of symmetric matrices. Moreover, it is assumed that $A_i,K_j,D_j$ are continuous in $\epsilon$ uniformly over bounded $U$ and that they are $C^s$ in $U$ uniformly in $\epsilon$ for $s-1\geq s_0:=\frac{N+M}{2}$. The initial value $U_0$ is assumed to be an element of $H^s(\R^N\times\mathbb{T}^M)$. For the reduced equation, we introduce the $C_0$-semigroup $(S(\tau))_{\tau\geq0}$ generated by the system
 \[
    \partial_{\tau}v + \sum_{j=1}^M K_j\partial_{y_j}v=0.
 \]
 Moreover, we define the inner product
 \begin{align*}
    &\langle u,v\rangle_{L_{2}([0,\infty),L_2(\R^N\times\mathbb{T}^M))} :=\\ &\qquad\qquad \lim_{\tau\to\infty}\frac{1}{\tau}\int_0^{\tau}\int_{\R^N\times\mathbb{T}^M} u(s,x,y)v(s,x,y)\,d(x,y)\,ds
 \end{align*}
 $L_{\infty}([0,\infty),L_2(\R^N\times\mathbb{T}^M))$ as well as the orthogonal projection $P$ with respect to this inner product to the closed subspace generated by expressions of the form $S(\tau)g(x,y)$ for $g\in H^{s-1}(\R^N\times\mathbb{T}^M)$. With this notation at hand, we may define the reduced equation given by
 \begin{align}\label{Eq:SymHypReduced}\begin{aligned}
     &\partial_t V^0 + P\bigg([\partial_U A^0(0)V^0]\partial_{\tau}V^0+\sum_{i=1}^N A_i(V^0,0)\partial_{x_i}V^0\\
     &\qquad\qquad+\sum_{j=1}^M D_j(V^0,0)\partial_{y_j} V^0- F(V^0,0)\bigg)=0,\\
    & V^0(0,x,\tau,y) = S(\tau) U_0(x,y).
 \end{aligned}\end{align}
 This is now a reduced equation in the following sense (see \cite[Theorem 2.3]{schochet1994fast}: Let $V^0$ be the solution of \eqref{Eq:SymHypReduced} and suppose that there is a $T>0$ such that the mapping
 \[
    \tilde{V}^0:t\mapsto V^0(t,\cdot_x,t/\epsilon,\cdot_y)
 \]
 is an element of $L_{\infty}([0,T];H^s(\mathbb{R}^N\times\mathbb{T}^M))$ for small $\epsilon>0$. Then, there are $\epsilon_0,C>0$ such that for all $\epsilon\in(0,\epsilon_0)$ and all $t_0\in[0,T]$ the solution $U^{\epsilon}$ of \eqref{Eq:SymHypSystem} exists and it holds that
 \[
    \|U^{\epsilon}-\tilde{V}^0\|_{C([0,t_0],H^{s-1}(\mathbb{R}^n\times\mathbb{T}^M))}\leq C t_0.
 \]
 This result can be improved under further assumptions: The crucial assumption is that the operator that maps $(\tau,x,y)\to S(\tau)U_0(x,y)$ to
 \begin{align*}
   & \left(\partial_{\tau}+\sum_{j=1}^M K_j\partial_{y_j}\right)^{-1} (1-P)\bigg(\partial_U A^0(0)V^0]\partial_{\tau}V^0+\\
    &\sum_{i=1}^N A_i(V^0,0)\partial_{x_i}V^0+\sum_{j=1}^M D_j(V^0,0)\partial_{y_j} V^0- F(V^0,0)\bigg)\bigg\vert_{\tau=t/\epsilon}
 \end{align*}
 for all $t\in[0,T]$ has to be bounded from $H^s(\R^n\times\mathbb{T}^M)$ to $H^{s-1-p}(\R^n\times\mathbb{T}^M)$ for some $p\in\R$. Now, if the other assumptions are satisfied even for $s-1\geq s_0+\max\{0,p\}$ and if $A_i,D_j,F$ are continuously differentiable in $\epsilon$, then it even holds that 
 \[
    \|U^{\epsilon}-\tilde{V}^0\|_{C([0,t_0],H^{s-1}(\mathbb{R}^n\times\mathbb{T}^M))}\leq C t_0\epsilon,
 \]
 see \cite[Corollary 2.4]{schochet1994fast}. We refer to \cite[Lemma 2.5]{schochet1994fast} for an analysis of these additional assumptions.\\
 Already in \cite{klainerman1981singular,schochet1994fast}, these techniques for singular limits of hyperbolic systems have been used to derive rigorous results for incompressible limits of fluid equations. They have been further refined in \cite{embid1996averaging}
 in which the so-called {\em fast wave averaging} has been developed. This fast wave averaging has been used in the literature to study the limiting behaviour for geophysical flows with different combinations of low or finite Froude and Rossby numbers, see \cite{embid1998low,wingate2011low}.
\end{itemize}

\subsection{Homogenization of Spatially Periodic Structures} \label{Sec:Homogenization_Spatially_Periodic}
A standard homogenization problem which is introduced in many textbooks is the one for elliptic partial differential equations of the form
\begin{align}\label{Eq:Homogenization:FastCoefficients}
     -\nabla\cdot(A^{\epsilon}\nabla u^{\epsilon})=f
    \end{align}
on a domain $\mathcal{O}\subset\R^n$ with Dirichlet boundary conditions. Here, $A_{\epsilon}=A(y/\epsilon)$ for some $1$-periodic, $\R^{n\times n}$-valued function $A$ such that \eqref{Eq:Homogenization:FastCoefficients} is uniformly elliptic. $f\colon\mathcal{O}\to\R$ is a given forcing term. For small values of $\epsilon>0$, the coefficients coming from $A_{\epsilon}$ have strong spatial oscillations. The aim is to reduce \eqref{Eq:Homogenization:FastCoefficients} to an equation of similar type, but without the oscillating coefficients. More precisely, it will be reduced to 
allows one to derive a reduced equation for $u_0$, namely
    \begin{align*}
     -\sum_{i,j=1}^n \overline{a}_{i,j} \partial_{y_i}\partial_{y_j}u_0=f,
    \end{align*}
 where the $\overline{a}_{i,j}$ are just scalars and where $u_0$ nonetheless is a good approximation of $u^{\epsilon}$ in \eqref{Eq:Homogenization:FastCoefficients}.
\smallskip

 {\bf{Outline of the Reduction Procedure:}}  A popular approach is to make the ansatz
 \begin{align*}
 u^{\epsilon}(y)&=u_0(\tfrac{y}{\epsilon},y)+\epsilon u_1(\tfrac{y}{\epsilon},y)+\epsilon^2u_2(\tfrac{y}{\epsilon},y)+\ldots\\
 &=u_0(x,y)+\epsilon u_1(x,y)+\epsilon^2u_2(x,y)+\ldots
 \end{align*}
so that one obtains an asymptotic expansion in $\epsilon$ together with a two-scale structure of the variables. Using this two-scale ansatz one can rewrite 
\[
 -\nabla\cdot(A^{\epsilon}\nabla)=:L=\frac{1}{\epsilon^2}L_1+\frac{1}{\epsilon}L_2+L_3
\]
where
\begin{align*}
 L_1w&:=-\sum_{i,j=1}^n\frac{\partial}{\partial x_i}a_{i,j}(x)\frac{\partial}{\partial x_j}w,\\
 L_2w&:=-\sum_{i,j=1}^n\frac{\partial}{\partial x_i}a_{i,j}(x)\frac{\partial}{\partial y_j}w+\frac{\partial}{\partial y_i}a_{i,j}(x)\frac{\partial}{\partial x_j}w,\\
 L_3w&:=-\sum_{i,j=1}^n\frac{\partial}{\partial y_i}a_{i,j}(x)\frac{\partial}{\partial y_j}w,
\end{align*}
with $A(x)=(a_{i,j}(x))_{1\leq i,j,\leq n}$.
Now one plugs the expansion for $u^{\epsilon}$ into \eqref{Eq:Homogenization:FastCoefficients} and separates the different orders in $\epsilon$. This yields the equations
\begin{align*}
 0&=L_1u_0,\\
 0&=L_1u_1+L_2u_0,\\
 f&=L_1u_2+L_2 u_1+ L_3 u_0.
\end{align*}
One obtains that the solution of $L_1u_0=0$ only depends on $y$ and is constant in $x$. Introducing the $1$-periodic solution $\chi_i$  of the cell problem
\begin{align*}
 L_1\chi_i=-\sum_{j=1}^n \frac{\partial}{\partial x_j}a_{i,j}(x)
\end{align*}
allows one to derive a reduced equation for $u_0$, namely
    \begin{align}\label{Eq:Homogenization:Reduced}
     -\sum_{i,j=1}^n \overline{a}_{i,j} \partial_{y_i}\partial_{y_j}u_0=f
    \end{align}
with Dirichlet boundary conditions, where
\[
 \overline{a}_{i,j}:=\int_{\mathbb{T}^n}a_{i,j}(x)-\sum_{k=1}^n a_{i,k}(x)\frac{\partial\chi_i}{\partial x_k}\,dx,
\]
i.e. the $\overline{a}_{i,j}$ are just constants. Finally, one can show that $u^{\epsilon}\to u_0$ in the $L_2(\mathcal{O})$-sense as $\epsilon\to0$, so that one may reduce \eqref{Eq:Homogenization:FastCoefficients} to \eqref{Eq:Homogenization:Reduced}.
 
\smallskip

 {\bf{Example:}}  One can use homogenization to derive Darcy's law, which models the laminar flow of a fluid in porous media, from Navier-Stokes type equations. There exist different variants of the reduction procedure. We refer to \cite{Lu_2020} and references therein for an overview of them. Here, we briefly describe the reduction in \cite{Mikelic_1991}. Therein, the incompressible Navier-Stokes equation
\begin{align}
    \begin{aligned}\label{Eq:Homogenization:Navier-Stokes}
    \partial_t v^{\epsilon} +(v^\epsilon\cdot\nabla) v^{\epsilon}+\nabla p^{\epsilon}-\mu\Delta v^{\epsilon} &= f^{\epsilon}\qquad\text{in } \Omega_{\epsilon}\times(0,T),\\
    \operatorname{div} v^{\epsilon}&=0\qquad\text{in } \Omega_{\epsilon}\times(0,T),\\
    v^{\epsilon}&=h\qquad\text{on } \Gamma\times(0,T),\\
    v^{\epsilon}&=0\qquad\text{on } S_{\epsilon}\times(0,T),\\
    v^{\epsilon}(x,0)&=v_0^{\epsilon}(x)\qquad\text{in } \Omega_{\epsilon}
    \end{aligned}
\end{align}
is studied. Here, the periodicity comes from the special shape of the domain $\Omega_{\epsilon}$ which is given as a smooth domain with periodically distributed holes. More precisely, one starts with a smooth domain $\Omega\subset\R^n$ with $\Gamma:=\partial\Omega$ and a second smooth domain $\mathcal{D}\subset[0,1]^n$. Now, one defines the set $T_{\epsilon}:=\{k\in\Z^n: \epsilon k+[0,\epsilon]^n\subset\Omega\}$ as well as
\[
    \mathcal{D}_{\epsilon}:=\bigcup_{k\in T_{\epsilon}} \epsilon(k+\mathcal{D}),\quad S_{\epsilon}:=\partial\mathcal{D}_{\epsilon},\quad \Omega_{\epsilon}:=\Omega\setminus\overline{\mathcal{D}}_{\epsilon},
\]
where $\overline{\mathcal{D}}_{\epsilon}$ denotes the closure of $\mathcal{D}_{\epsilon}$. The given data is assumed to satisfy
\begin{align*}
    h&\in C([0,T];W^{1}_{\infty}(\Omega;\R^n)),\quad \partial_th = C([0,T];W^1_6(\Omega;\R^n)),\quad\operatorname{div} h = 0,\\
   &\qquad \epsilon \| v_0^{\epsilon} \|_{L_2(\Omega_{\epsilon};\R^n)} = \mathcal{O}(1),\quad \operatorname{div} v_0^{\epsilon} = 0\quad\text{in }\Omega_{\epsilon},\\
  &\qquad f^{\epsilon}\in L_2((0;T)\times\Omega_{\epsilon};\R^n),\quad \epsilon^{2}f^{\epsilon}\to f\quad\text{as }\epsilon\to 0,
\end{align*}
for a certain $f \in L_2((0;T)\times\Omega_{\epsilon};\R^n)$. Note that $h$, which defines the flux through the outer boundary $\Gamma$, is chosen such that 
\[
    \int_{\Gamma} \langle h(x),\nu(x) \rangle\,\txtd S(x)=0,
\]
where $\nu$ denotes the unit outer normal on $\Gamma$ and $S$ denotes the surface measure.\\
For the reduction to Darcy's law, the solutions of a suitable cell problem play again an important role. This time, it takes the form
\begin{align}
    \begin{aligned}\label{Eq:Homogenization:DarcyCellProblem}
        -\nabla \pi^j + \Delta w^j + e^j& = 0 \qquad\text{in }[0,1]^n\setminus \mathcal{D},\\
          \operatorname{div} w^j&=0\qquad\text{in }[0,1]^n\setminus \mathcal{D},\\
          w^j&=0\qquad\text{on } \partial\mathcal{D},
    \end{aligned}
\end{align}
with periodic boundary conditions on $\partial[0,1]^n$, where $e^j$ denotes the $j$-th canonical unit normal vector in $\R^n$. From the unique solution $(w^j,\pi^j)$ (up to constants in $\pi$) of \eqref{Eq:Homogenization:DarcyCellProblem} one can construct the so-called permeability tensor $K=(K_{j,k})_{1\leq j,k\leq n}$ by means of
\[
    K_{j,k}=\int_{[0,1]^n\setminus \mathcal{O}} (w^j)_k(x)\,\txtd x.
\]
This permeability tensor now appears in the homogenized equation which is given by
\begin{align}
    \begin{aligned}\label{Eq:Homogenization:Darcy}
           \operatorname{div} v&=0\qquad\text{in } L_2(\Omega\times(0,T)),\\
          v=&=\frac{K}{\mu}(f-\nabla p)\qquad\text{in } L_2(\Omega\times(0,T)),\\
          v\cdot\nu&=h\cdot\nu\qquad\text{in } L_2((0,T);H^{1/2}(\Gamma)).
    \end{aligned}
\end{align}
Solutions $(v^{\epsilon},p^{\epsilon})$ of \eqref{Eq:Homogenization:Navier-Stokes} are defined on $\Omega_{\epsilon}$. However, the aim of the homogenization procedure is to average out the effect of the holes so that \eqref{Eq:Homogenization:Darcy} and its solutions $(v,p)$ should be defined on $\Omega$. Therefore, if one wants to compare the approximation $(v,p)$ to the solution of the original problem $(v^{\epsilon},p^{\epsilon})$, one should either restrict $(v,p)$ to $\Omega_{\epsilon}$ or extend $(v^{\epsilon},p^{\epsilon})$ in a suitable way to $\Omega$. The latter approach, which has the advantage that the approximation result can be formulated in terms of convergence in a fixed space with a fixed topology, was chosen in \cite{Mikelic_1991}. While $v^{\epsilon}$ can just be extended by $0$, the adjoint of a certain restriction operator is used to define the extension $\tilde{P}^{\epsilon}$ of a time integrated version $P^{\epsilon}$ of the pressure $p^{\epsilon}$. We refer to \cite[Section 3]{Mikelic_1991} for the details of the construction. By \cite[Theorem 4.1]{Mikelic_1991} the extensions $v^{\epsilon}$ and $\tilde{P}^{\epsilon}$ now satisfy
\begin{align*}
    v^{\epsilon}\to v\quad\text{weakly in } L_2 ((0,T)\times\Omega;\R^n),\\
    \epsilon^{2} \partial_t \tilde{P}^{\epsilon} \to p\quad\text{weakly in } W^{-1}_2 ((0,T);L_{\beta}^0(\Omega)),
\end{align*}
where $1<\beta < \frac{n}{n-1}$ and where $L_{\beta}^0(\Omega):=\{f\in L_{\beta}(\Omega):\int_{\Omega} f(x)\,\txtd x=0\}$.

\smallskip

 {\bf{Comments:}}  \begin{itemize}
                 \item An important concept for homogenization techniques is the so-called two-scale convergence, as it is used in many proofs related to homogenization. Let $\mathcal{O}\subset\R^n$ be a domain and $\mathbb{T}^n$ the $n$-dimensional torus. Then we say that a family of functions $(u^{\epsilon})_{\epsilon\in(0,1)}\subset L_2(\mathcal{O})$ two-scale converges to $u^0\in L_2(\mathcal{O}\times \mathbb{T}^n)$ if for all $\phi\in L_2(\mathcal{O};C(\mathbb{T}^n))$ it holds that
\[
	\lim_{\epsilon\to 0} \int_{\mathcal{O}} u^{\epsilon}(x)\phi(x,\tfrac{x}{\epsilon})\,dx = \int_{\mathcal{O}}\int_{\mathbb{T}^n} u^0(x,y)\phi(x,y)\,dy\,dx.
\]
\item Homogenization techniques also exist for parabolic partial differential equations. They are usually applied to equations of the form
\begin{align}
	\begin{aligned}\label{Eq:Homogenization:Parabolic}
		\partial_t u^{\epsilon}(t,x)&= \frac{1}{\epsilon} \langle b\big(\tfrac{x}{\epsilon}\big),\nabla\rangle u^{\epsilon}(t,x) + D\Delta u^{\epsilon}(t,x)\quad((t,x)\in \R_+\times\R^n),\\
		u^{\epsilon}(0,x)&=u_0(x)\quad(x\in\R^n),
	\end{aligned}
\end{align}
where $D>0$ and $b\colon\R^n\to\R^n$ is smooth and $1$-periodic in all directions. Such equations are for example treated in \cite[Chapter 13 \& Chapter 20]{Pavliotis_Stuart_2008}.
\end{itemize}

Yet another important variant of homogenization in the context of climate dynamics are techniques reducing fast chaotic degrees of freedom to a lower-dimensional stochastic input, resulting in a reduced model.

\subsection{Stochastic modelling of Deterministic Systems via Homogenization}
\label{Sec:Homog}

It is a common idea in climate dynamics that variability corresponding to chaotic parts of the dynamics which happen on a fast time scale can effectively be approximated by stochastic terms. This idea goes back to Hasselmann \cite{Hasselmann_1976} and has since been used in many climate models, with significant contributions from Majda and co-workers \cite{Majda_Timofeyev_VandenEijnden_2001,Majda_Timofeyev_VandenEijnden_2003,Majda_Franzke_Crommelin_2009}. From a mathematical perspective, not all the reductions from multi-scale ordinary differential equations to stochastic differential equations are rigorously justified. Nonetheless, there are mathematical techniques which can make the reductions rigorous. These techniques are usually also called homogenization, but are different from the ones in Section~\ref{Sec:Homogenization_Spatially_Periodic}. The starting point of such homogenization techniques are fast-slow systems of ordinary differential equations of the form

\begin{align}
    \begin{aligned}\label{Eq:Fast-Slow_Multiscale}
        \frac{du^{\epsilon}}{dt}&= \frac{1}{\epsilon}f_0(u^{\epsilon},v^{\epsilon})+ f_1(u^{\epsilon},v^{\epsilon}),\\
        \frac{dv^{\epsilon}}{dt}&= \frac{1}{\epsilon^2}g_0(u^{\epsilon},v^{\epsilon})+\frac{1}{\epsilon}g_1(u^{\epsilon},v^{\epsilon})+g_2(u^{\epsilon},v^{\epsilon}),
    \end{aligned}
\end{align}
with nonlinearities $f_0,f_1,g_0,g_1,g_2$ and $u^{\epsilon}$ representing slow variables forced by fast variables $v^{\epsilon}$. The aims is to find reduced equations are stochastic differential equations of the form
\begin{align}\label{Eq:Homogenized_SDE}
   dU(t)=\overline{g}(U(t))\,dt+\overline{\Sigma}(U(t))\,dW(t),
\end{align}
with nonlinearities $\overline{g}$ and $\overline{\sigma}$, and a certain noise term $W$ which is a Brownian motion in many but not in all cases. It is important to note that unlike in \eqref{Eq:Averaging:FastSlow_Generalized} the equations we start with are deterministic, i.e., we do not assume that $\sigma_1\sigma_1^T$ is uniformly positive definite. Rigorous reduction results from \eqref{Eq:Fast-Slow_Multiscale} to \eqref{Eq:Homogenized_SDE} are much more recent than the ones for \eqref{Eq:Averaging:FastSlow_Generalized}. We refer to \cite{Engel_Gkogkas_Kuehn_2021,Gottwald_Melbourne_2013,Kelly_Melbourne_2017,Melbourne_Stuart_2011}.\smallskip

{\bf{Outline of the Reduction Procedure:}} As already mentioned above, a system with multiple time scales of the form \eqref{Eq:Fast-Slow_Multiscale} is reduced to a stochastic differential equation of the form \eqref{Eq:Homogenized_SDE} with different choices of $\overline{g}$, $\overline{\sigma}$ and $W$, depending on the assumptions. In fact, one may even have to use different kinds of stochastic integrals.

This reduction is shown in a simpler setting in \cite{Melbourne_Stuart_2011} starting with
\begin{align}
    \begin{aligned}\label{Eq:Fast-Slow_Multiscale_MS2011}
        \frac{du^{\epsilon}}{dt}&= \frac{1}{\epsilon}f_0(v^{\epsilon})+f_1(u^{\epsilon},v^{\epsilon}),\\
        \frac{dv^{\epsilon}}{dt}&= \frac{1}{\epsilon^2}g_0(v^{\epsilon}),
    \end{aligned}
\end{align}
for $u\in\R^d$, with initial condition $(u^{\epsilon}(0),v^{\epsilon}(0))=(\xi,\eta)$,  there is an invariant measure $\mu$ for the fast dynamics $v^{\epsilon}$ that is assumed uniformly or nonuniformly hyperbolic with fast decay of correlations and $f_0$ has zero mean with respect to $\mu$. It is then possible to show \cite[Theorem 1.1]{Melbourne_Stuart_2011} that under reasonable regularity and dynamical assumptions there is weak convergence of solutions $u^{\epsilon}\rightarrow_w U$ for $t>0$ as $\epsilon\rightarrow 0$, where $U$ satisfies an SDE on $\R^d$ given by
$$
U(t)=\xi+\int_0^{t} F(U(s))\,ds+\sqrt{\Sigma}W(t)
$$
where $W$ is Brownian motion, $\Sigma$ a covariance matrix and $F$ is obtained by averaging $f$ with respect to $\mu$. This relies on assumptions on the fast dynamics that ensure there is a Weak Invariance Principle (WIP) of the form 
$$
\frac{1}{\sqrt{n}} \int_0^{nt} f_0(v^{1}(\tau))\,d\tau \rightarrow_w \sqrt{\Sigma} W(t)
$$
in $C^{0}([0,T],\R^d)$ as $n\rightarrow\infty$ as well as a large deviation assumption on the fast dynamics. 

This is generalized in \cite{Gottwald_Melbourne_2013,Kelly_Melbourne_2017} to cases that include not only additive but more general (such as multiplicative) noise, where $f_0$ depends not just on $v^{\epsilon}$ but also on $u^{\epsilon}$, including cases (with slow decay of correlations) where the limiting SDE is driven by a stable L\'{e}vy process rather than a Wiener process. Interesting further work in this direction includes \cite{Wouters_Gottwald_2019} who use Edgeworth expansions to find corrections to asymptotic homogenization results for small $\epsilon$.\smallskip

{\bf{Comments:}}  \begin{itemize}
\item Although, since the work of Hasselmann  \cite{Hasselmann_1976}, homogenization has arguably been one of the most applied methods in climate dynamics, it is also one of the methods that is most difficult to justify in a rigorous manner. Results so far tend to assume a skew product structure, quite strong assumptions on the forcing chaos and also a timescale separation close to an asymptotic limit. All of these assumptions cannot be verified except in some fairly limited settings, however the utility of the stochastic modelling approach goes well beyond what can be rigorously justified.
\item The methods discussed in Section~\ref{sec:Slow} as well as the methods discussed in this section do rely on scale separation. Yet, there are interesting problems, where scale separation breaks down. For example, in the context of time scales, classical non-autonomous dynamical systems theory~\cite{KloedenRasmussen} provides tools to formulate and understand relevant notions of attraction and attractor in the non-autonomous context, in particular pullback notions of attraction. In climate dynamics and other applications, this has gained recent attention under the theme of rate-induced instability effects~\cite{Ashwin_2012,KuehnLongo,alkhayuon2019basin,Wieczorek_2021} which presents some methods for understanding the location where scale separation breaks, and what happens at such rate-induced tipping points.
\end{itemize}

Related to averaging and homogenization is another classical idea from statistical physics, where a reduction is performed using the concept of averaging in a slightly different variant.

\subsection{Mori-Zwanzig Reduction}
\label{Sec:Mori}

Another reduction method that can be applied to remove fast degrees of freedom from an evolution equations is the projection method originally developed by Mori and Zwanzig \cite{Mori_1965,Zwanzig_1961,Zwanzig_1973} and for example discussed in \cite{Chorin_2002}. This method has recently found use in several climate applications \cite{Wouters_Lucarini_2013,Falkena_etal_2019}. It is unusual in relation to the other techniques discussed in this review in that it includes a memory term.

{\bf Outline of Reduction Procedure:}
If we consider an evolution equation of the form
$$
\frac{d}{dt}\varphi= R(\varphi),~~\varphi(0)=x
$$
that generates a flow $\varphi(x,t)$, for the sake of simplicity we assume $\varphi$ (and therefore $x$) are in $\R^n$. Suppose we would like to predict the evolution of some observable $u(x,t)=g(\varphi(x,t))$. Then by considering the linear PDE
\begin{equation}
\frac{\partial}{\partial t} u(x,t) = \cL u(x,t),~u(x,0)=g(x)
\label{eq:Liouville}
\end{equation}
where $[\cL u](x)=\sum_{i} R_i(x)\partial_{x_i}u(x)$, $\cL$ is the Liouville operator, and the solution of (\ref{eq:Liouville}) can be written 
$$
u(x,t)=[\exp(t\cL) g](x).
$$
We aim to understand the evolution of $g(\varphi(x,t))$ in terms of the evolution under (\ref{eq:Liouville}). Consider now a projection $P$ of $\R^n$ onto an $m$-dimensional subspace, with $m<n$. We wish to understand the dynamics of the observable using projection by $P$ onto a number of resolved variables $\hat{\varphi}:=P\varphi$, noting the unresolved variables can be written $\bar{\varphi}=(1-Q)\varphi$ (various choices for $P$ are discussed in \cite{Chorin_2002}).

For example, if we consider $g(x)=x_i$ as projection onto a single component then it is possible to show (for example, \cite{Chorin_2002,Falkena_etal_2019}) that it evolves according to a generalized Langevin-type equation of the form
\begin{equation}
\frac{\partial}{\partial t} \varphi_i(x,t) = R_i(\hat{\varphi}(x,t))+\int_{0}^t K_i(\hat{\varphi}(x,t-s),s)\,ds+F_i(x,t)
\label{eq:MZ}
\end{equation}
where 
$$
F_i(x,t)= [\exp(tQ\cL) Q\cL x]_i,~~K_i(\hat{x},t)=[P\cL F]_i(x,t).
$$
Note that for the resolved variables, the first term on the right hand side of (\ref{eq:MZ}) corresponds to a Markovian term that only depends on the current values of the resolved variables, the second term embodies a memory kernel while the third and final term encapsulates the influence of the unresolved variables. Depending on the situation, one may end up with an ODE reduction (if the last two terms in (\ref{eq:MZ}) are small) a DDE reduction (if the last term in (\ref{eq:MZ})  is small) or an SDE reduction (in cases where the middle term in (\ref{eq:MZ}) is small and the last term can be understood as a stochastic term).

{\bf Examples:}
The Mori-Zwanzig reduction method can be a useful method to derive a delay differential equation from a more complex PDE with hyperbolic terms. For example, \cite{Falkena_etal_2019} use this formalism to justify a rational reduction of a PDE model of the El Ni\~{n}o Southern Oscillation (ENSO) to a delay differential equation, while \cite{Falkena_etal_2021} use the method to justify a reduction of a PDE model of Atlantic Multidecadal Variability (AMV) to a delay differential equation.

\section{Moment Closure Methods}
\label{sec:mc}

The manifold reductions as well as techniques related to averaging and homogenization described in Sections~\ref{sec:linear}-\ref{Sec:AvHom} directly work with the equations describing a certain system. However, if a system is very complex one may prefer to only work with a set of observables of the original system. In complex nonlinear systems it is common that such a set of observables depends on other observables again. If one wants to describe all the observables through a closed set of differential equations, then one might obtain infinitely many observables with infinitely many equations. 

The following reduction method, the so-called moment closure methods, treat exactly such a kind of problems. Moment closure methods have been used in many different scientific disciplines, including climate \cite{franzke2005low-order}. They can be applied if one works with an abstract evolution equation of the form
\begin{align}\label{Eq:MomentClosure:Equation}
	\partial_t u = N(u).
\end{align}
This equation is to be understood in a formal way. For example, we also allow the case in which $N$ contains noise terms. Frequently, the method is applied to stochastic differential equations~\cite{Socha} of the form
\begin{align}\label{Eq:MomentClosure:SDE}
	dx_t=f(x_t)\,dt+F(x_t)\,dW_t.
\end{align}
But also other very important scenarios are possible such as kinetic PDEs~\cite{Levermore} of the form
\begin{align}\label{Eq:MomentClosure:Kinetic}
 \partial_t\varrho+v\cdot\nabla_x\varrho=Q(\varrho),
\end{align}
where $\rho=\rho(t,x,v)$ can be interpreted - if normalized - as the probability density of a single particle being located at $x$ and having velocity $v$ at time $t$. The so-called collision operator $Q$ usually only depends on the velocity $v$. Moment closure methods are also frequently used in network dynamics~\cite{GrossSayama}. We also refer to the references in the review~ \cite{Kuehn_2016}.  

\smallskip

 {\bf{Idea:}} An equation like \eqref{Eq:MomentClosure:Equation} may be too complex to be studied analytically or even numerically. However, one can still try to describe certain quantities extracted from the full system, the so-called moments. Usually, the moments are scalar-valued and have a certian hierarchy in the sense that there is a natural ordering which allows one the speak of higher order and lower order moments. One may use \eqref{Eq:MomentClosure:Equation} to derive new equations describing the dynamics of the moments. These new equations are referred to as moment equations. Technically, this step is not a reduction procedure in the sense that the moment equations do not approximate the original system. Instead, only certain aspects of the original system are considered, which in exchange are described at full resolution and not only approximated. The approximation procedure is contained in the next step: Oftentimes, the moment equations will be a fully coupled system of infinitely many differential equations and such a system may be as complicated as the original equation \eqref{Eq:MomentClosure:Equation}. In order to reduce the complexity, the moment closure is performed, i.e., the higher order moments are assumed to depend in a certain way on the lower order moments so that one can derive a closed system with finitely many variables. If the moment closure is performed in the right way, then the reduced system should be much simpler to study while still being a good approximation of the full system of moment equations. While in the abstract framework one might think that it is hard to find a good moment closure, many straightforward choices in applications work surprisingly well.

\smallskip

 {\bf{Outline of the Reduction Procedure:}} 
  The procedure usually consists of four steps 
    \begin{enumerate}
     \item One has to determine a set of moments $\{m_k:k\in\N\}$ one wants to consider. In applications, moments are often of the following form: One starts with a suitable solution concept for \eqref{Eq:MomentClosure:Equation} and a solution $u\colon\R_+\to X$ with some state space $X$. Moreover, one has a set of mappings $\{M_k:k\in\N\}$ with $M_k\colon X\to X$ as well as a mapping $\langle\cdot\rangle\colon X\to\mathbb{K}$, where $\mathbb{K}\in\{\R,\C\}$, which can be interpreted as a notion of average. One may think of $\langle\cdot\rangle$ as being a tool to remove complexity from the system, for example in the form of randomness or spatial dependence. The mappings $M_k$ in turn are used to still preserve the information on certain aspects of the system. Now, one defines the moments $m_k(t):=\langle M_k(u(t)) \rangle$. The precise choice of $\{M_k:k\in\N\}$ and the mapping $\langle\cdot\rangle$ depends a lot on the specific problem. The main requirement is that it should be possible to carry out the second step, the derivation of a system of moment equations, which we briefly explain below.\\
     One of the most common scenarios, which also motivates the terminology in the context of moment closure methods, is the case of a stochastic differential equation \eqref{Eq:MomentClosure:SDE}, where one may oftentimes just use the usual moment of a random variable. The solution to \eqref{Eq:MomentClosure:SDE} will be a stochastic process, and if it takes values in $X=\bigcap_{p>0} L_p(\Omega,\mathcal{F},\mathbb{P})$, where $(\Omega,\mathcal{F},\mathbb{P})$ denotes the underlying probability space, then one may define 
    \[
    	 M_k  \colon X\to X, \; u\mapsto u^k\quad\text{and} \quad \langle\cdot\rangle \colon X\to\R,\;u\mapsto \mathbb{E}[u].
    \]
    Therefore, $m_k$ is just the $k$-th moment of the solution process.
     \item One has to derive the moment equations. The details of this step depend a lot on the choice of moments in the first step. But generally, the aim is to derive an infinite system of differential equations for $\{m_k:k\in\N\}$, i.e. a system of the form
          \begin{align}
        \begin{aligned}\label{Eq:Moment_Closure:Moment_Equations}
      \dot{m}_1&=h_1(m_1,m_2,m_3,\ldots),\\
      \dot{m}_2&=h_2(m_1,m_2,m_3,\ldots),\\
      \dot{m}_3&=h_3(m_1,m_2,m_3,\ldots),\\
      \ldots&=\qquad\ldots,
      \end{aligned}
     \end{align}
     for certain $h_k\colon \mathbb{K}^{\N}\to \mathbb{K}$, $k\in\N$. Since we defined $m_k(t)=\langle M_k(u(t)) \rangle$, one might already guess that one would have to combine a chain rule together with the fact that $u$ solves \eqref{Eq:MomentClosure:Equation}. Indeed, this approach usually leads in the right direction, but in practice, this may cause difficulties. For example, $u$ might not be differentiable as in the case of a stochastic differential equation. Fortunately, one can still try to use It\^o's lemma in such a situation. Further below, we will carry out this step for the example of a stochastic version of the Stommel-Cessi model.\\
     Even though $h_k$ will not actually depend on all the moments in most applications, the system of moment equations can still be an infinite system, since the $h_k$ usually depend on higher order moments $m_{k+n_1},\ldots,m_{k+n_l}$ for some natural numbers $l,n_1,\ldots,n_l\in\N$. The problem of reducing this system to a closed finite systems of ordinary differential equations will be the main task of the third step.
     \item One has to close the system of moment equations in a suitable way. There are many ways to do this. In general, one looks for a mapping $H\colon\R^{N}\to\R^{\N}$ with suitable $N\in\N$ that determines the higher order moments through the lower order moments, i.e.
     $$H(m_1,\ldots,m_N)=(m_{N+1},m_{N+2},m_{N+3},\ldots).$$
     Naive approaches which are used quite often and which still work very well are $0=m_{N+1}=m_{N+2}=\ldots$ or more generally constant higher moments. One can also use certain assumptions coming from the application or the expected outcome. In the case of a stochastic differential equation, one can for example assume that one is looking for random variables with a certain distribution such as Gaussian random variables.
     \item Finally, one should justify that at least the lower order moments are approximated well by the closure procedure. However, since there are many different ways to carry out the first three steps, there is also no general rule for the last step. 
    \end{enumerate}
    
\smallskip

 {\bf{Example:}} As mentioned above, the reduction of stochastic differential equations is a standard application of moment closure methods. An example in climate dynamics is a stochastic version of the Stommel-Cessi model (see \cite{Cessi_1994,Stommel_1961}) given by
  \begin{align}\label{Eq:Stochastic_Stommel_Cessi}
     \txtd U_t = \mu - U_t [ 1 + 7.5 (1-U_t)^2]\,\txtd t + \sigma \,\txtd W_t.
  \end{align}
  Here, $U_t$ models the salinity difference between a higher and a lower latitude box in the North Atlantic, $\mu$ models the difference in freshwater flux, $(W_t)_{t\geq0}$ denotes a Brownian motion, and $\sigma$ is a small parameter which models the strength of the noise. This equation is derived from the Stommel-Cessi model by a reduction to the critical manifold and by adding the stochastic forcing, see for example \cite[Section 6.2.1]{Berglund_Gentz_2006} or \cite[Example 19.9.3]{Kuehn_2015}. In order to carry out a moment closure procedure, we define
    \[
    	 M_k  \colon X\to X, \; u\mapsto u^k\quad\text{and} \quad \langle\cdot\rangle \colon X\to\R,\;u\mapsto \mathbb{E}[u].
    \]
    with $X=\bigcap_{p>0} L_p(\Omega,\mathcal{F},\mathbb{P})$ as above. For the moment equations, we use It\^o's Lemma and obtain
    \begin{align*}
        \txtd U^k_t = &k U_t^{k-1} \txtd U_t +\frac{k(k-1)}{2}U_t^{k-2} \txtd [U_t] \\
        = & k \big[-7.5U_t^{k+2}+15U_t^{k+1}-8.5U_t^k+k\mu U_t^{k-1}\\
        & +\frac{(k-1)\sigma^2}{2}U_t^{k-2}\,\txtd t + \sigma U_t^{k-1} \,\txtd W_t \big]
    \end{align*}
    for $k\in\N$, $k\geq2$. Taking the expectation and differentiating with respect to $t$ yields
    \[
        \dot{m}_k(t)= \big[-7.5m_{k+2}(t)+15 m_{k+1}(t)-8.5m_k(t)+k\mu m_{k-1}(t)+\frac{(k-1)\sigma^2}{2}m_{k-2}(t)\big]
    \]
    Using the convention $m_{0}=1$ and $m_{-1}=0$, this equation is valid for all $k\in\N$ and yields the system of moment equations.
    
    \smallskip
    
  {\bf{Comments:}}
  \begin{itemize}
      \item Formulating abstract results on the validity of moment closure methods still seems to be an open problem. However, moment closure methods often work very well in practice, when they are benchmarked via numerical simulations.
      \item In some sense, we may view moment closure methods as related to other reduction methods. The idea to use set of basis functions from Galerkin or the idea to use a linear subspace from EOF could also be viewed as one choice of the space of observables. The idea to use averaged observables is connected to previously presented ideas on averaging and homogeneization, while selecting a suitable closure is connected to invariant manifolds~\cite{Kuehn_2016}. 
  \end{itemize}
  
Most of the reduction methods we discussed so far aim to reduce the complexity or dimension of a dynamical system which is given by an explicit equation. However, nowadays large amounts of data on the Earth's climate system is collected so that compressing and finding structures in large data sets can be considered as least as important for the understanding of the Earth's climate as the formal reduction of dynamical systems. In principle, the EOFs approach described in Section~\ref{Sec:EOFs} fulfils this purpose very well, but as a linear method, it may have difficulties in preserving nonlinear structures such as lower-dimensional manifolds. Diffusion maps, which we discuss in the next section, are designed for exactly that purpose and are just one, of the many, data-driven reduction methods that start to permeate many applications in climate science.

\section{Diffusion Maps}
\label{Sec:DiffusionMaps}

Diffusion maps are a relatively young and mainly data-driven method of dimensionality reduction, which has first been introduced in \cite{Coifman_Lafon_2006}. Formally, it can be applied whenever EOFs can be used. But while EOFs are a linear dimensionality reduction method, and as such have problems in finding lower-dimensional but nonlinear structures in the data set, diffusion maps specifically aim at finding low-dimensional structures in high-dimensional data and representations thereof in a low-dimensional Euclidean space.
   
\smallskip

 {\bf{Idea:}} On a set of data one defines a kernel which is then used to construct a Markov chain with the data set as a state space. The Markov chain can be thought of as a form of diffusion which takes the geometry of the data set into account. The hope is to get some insight on this geometry by studying the structure of the Markov chain. This is achieved by defining the so-called diffusion maps, which are constructed from the eigenvalues and eigenfunctions of the Markov chain. The diffusion maps embed the data set into Euclidean space such that the Euclidean distance of two embedded points is close to the distance of the original points in the geometry of the data set.
    
\smallskip

 {\bf{Outline of the Reduction Procedure:}} One starts with a measure space $(X,\mathscr{A},\mu)$ and a map $k\colon X\times X\to [0,\infty)$ which is symmetric in the sense that $k(x,y)=k(y,x)$ for $x,y\in X$. This map is thought to measure the distance between two points in the data set $X$: the larger $k(x,y)$ the closer are $x,y$ thought to be. Technically, one has to impose further assumptions on $k$ so that the objects we define in the following are well-defined. But since these assumptions are very mild and satisfied in most situations, we refrain from being precise here. Instead, we refer to \cite{Coifman_Lafon_2006} for more details on the construction.\\
 Having chosen $k$ one defines
 \[
 	d(x):=\int_{X}k(x,y)\,d\mu(y)
 \]
 which can be interpreted as a local measure of volume. This local measure of volume is used to renormalize $k$ and define
 \[
 	p(x,y):=\frac{k(x,y)}{d(x)}.
 \]
 The mapping $p$ is not symmetric anymore, but it satisfies $\int_X p(x,y)\,dy=1$. It can therefore be seen as the transition kernel of the Markov chain $((P^*)^n)_{n\in\N}$ given as the powers of operator
 \[
 	P^*\colon L_2(X,\mathscr{A},\mu)\to L_2(X,\mathscr{A},\mu),\;P^*f(y):=\int_X p(x,y) f(x)\,d\mu(x).
 \]
 This operator is the adjoint of the so-called diffusion operator
  \[
 	P\colon L_2(X,\mathscr{A},\mu)\to L_2(X,\mathscr{A},\mu),\;Pf(x):=\int_X p(x,y) f(y)\,d\mu(y).
 \]
 which leaves constant functions invariant due to $\int_X p(x,y)\,dy=1$. Under the right conditions on $k$ (again we refer to \cite{Coifman_Lafon_2006}) there is a unique stationary distribution of the Markov chain given by
 \[
 	\pi(x):=\frac{d(x)}{\int_X d(z)\,d\mu(z)}.
 \]
 Using this stationary distribution one can define the modified kernel
 \[
 	a(x,y):=\frac{\sqrt{\pi(x)}}{\sqrt{\pi(y)}}p(x,y)=\frac{k(x,y)}{\sqrt{\pi(x)}\sqrt{\pi(y)}}.
 \]
 This kernel is symmetric so that the corresponding operator
   \[
 	A\colon L_2(X,\mathscr{A},\mu)\to L_2(X,\mathscr{A},\mu),\;Af(x):=\int_X a(x,y) f(y)\,d\mu(y)
 \]
 is self-adjoint. If $\int_X\int_X a(x,y)^2\,d\mu(x)\,d\mu(y)<\infty$, then $A$ is even compact and there is a sequence of eigenvalues $1=\lambda_0\geq|\lambda_1|\geq|\lambda_2|\geq\ldots$ and orthonormal eigenfunctions $(\eta_l)_{l\in\N_0}$ in $L_2(X,\mathcal{A},\mu)$ such that
 \[
	a(x,y)=\sum_{l\in\N_0} \lambda_l\eta_l(x)\eta_l(y).
 \]
 	If we now define
 \[
 	\psi_l(x):=\frac{\eta_l(x)}{\sqrt{\pi(x)}},\quad\varphi_l(y):=\eta_l(y)\sqrt{\pi(y)}
 \]
 then we have the representation
  \[
	p(x,y)=\sum_{l\in\N_0} \lambda_l\psi_l(x)\varphi_l(y).
 \]
 Moreover, the $t$-th powers of $P$ are given as an integral operator with kernel
   \[
	p_t(x,y)=\sum_{l\in\N_0} \lambda_l^t\psi_l(x)\varphi_l(y).
 \]
 With all this notation at hand, we can define the diffusion distances $D_t$  by
\begin{align*}
	D_t(x,y):&=\left(\int_X\frac{(p_t(x,u)-p_t(y,u))^2}{\pi(u)}\,d\mu(u)\right)^{1/2}\\
	&=\left(\sum_{l\in\N_0} \lambda_j^{2t}(\psi_l(x)-\psi_l(y))^2\right)^{1/2}.
\end{align*}
For fixed $t\geq0$ this defines a metric on the data set $X$. This expression is small if there is a large number of short paths between $x$ and $y$, i.e. if there is a high probability that the Markov chain ends up in $y$ after $t$ steps if it starts in $x$. In this sense, the diffusion distance contains a lot of information on the internal structure of the data set $X$.\\
The reduction is now carried via the so-called diffusion maps given by
\[
 \Psi_t\colon X\to \R^N,\,~\Psi_t(x)=\begin{pmatrix} \lambda_0^t\psi_0(x)\\ \vdots\\ \lambda_{N}^t\psi_N(x)\end{pmatrix}
\]
up to a desired dimension $N$. In \cite{Coifman_Lafon_2006} the authors include a certain rule of the dimensionality that should be considered. They take $N=s(\delta,t):=\max\{l\in\N:|\lambda_l|^t>\delta|\lambda_1|^t\}$. This way, the diffusion maps preserve the diffusion distance to a certain extent in the sense that
\[
	\| \Psi_t(x)-\Psi_t(y) \|_{2}=\left(\sum_{l=0}^{s(\delta,t)} \lambda_j^{2t}(\psi_l(x)-\psi_l(y))^2\right)^{1/2}:=D_{t,\delta}(x,y),
\]
where $\|\cdot\|_2$ denotes the usual Euclidean distance. Therefore, the diffusion maps give a representation of the data set $X$ in an $N$-dimensional Euclidean space, such that the internal notion of distance in $X$ is quantitatively preserved in $\R^N$ with the Euclidean distance up to a small error.
   
\smallskip

{\bf{Examples of Applications:}} Since diffusion maps are mainly a data-driven method, the results one obtains from the reduction procedure are probably most interesting when it is applied to real world data. Instead of explicitly carrying out such an example, let us briefly explain a few common scenarios. A standard situation is if $X=\{x_1,\ldots, x_K\}\subset\R^d$ is a finite set of points in $\R^d$. Then the $\sigma$-field would just be the power set $\mathscr{P}(X)$ of $X$ and $\mu$ could for example be the counting measure $\zeta$. In this case, the integrals just turn into sums, i.e.
\[
 \int_X f(x)\,d\mu(x)=\sum_{i=1}^Kf(x_i).
\]
A common choice for a kernel in this situation would be $$k_{\epsilon}(x,y)=\exp\bigg(\frac{\|x-y\|^2}{2\epsilon}\bigg)$$ for a parameter $\epsilon>0$. With these choices one could technically apply a diffusion map reduction whenever an EOF-based reduction is also possible. As for EOFs, one can also use diffusion maps to study dynamical systems by applying them to simulated data from the system. As an example, we mention \cite{Coifman_et_al_2016}.

\smallskip

An interesting use case is if $X:=M$ is a compact submanifold $\R^d$ (possibly with boundary) with its Riemannian measure $\nu$. In practice, $X$ might only consist of a finite number of sample points with a non-uniform distribution on $M$. Although we focus first on the continuous case, i.e. $X=M$, instead of the discrete case, we adopt the idea of a non-uniform distribution and introduce a density $q\colon M\to (0,\infty)$ on $M$. Then we take $d\mu=q\,d\nu$. Obviously, the outcome of the diffusion map reduction can depend on the density $q$ as well as the geometry of $M$. By choosing the right kernels, one can control the effect of the density. Depending on what is thought to be more important in a particular application, one can choose a corresponding kernel. More precisely, the kernels are constructed as follows: One starts with rotation-invariant reference kernels 
\[
 k_{\epsilon}(x,y):=h(\|x-y\|/\epsilon),
\]
where $h\colon[0,\infty)\to (0,\infty)$ is an exponentially decaying function, for example $h(x)=e^{-x^2/2}$ as above. An approximation of the density is defined by
\[
 q_{\epsilon}(x):=\int_X k_{\epsilon}(x,y) q(y)\,dy.
\]
Now, a new family of kernels is constructed by
\[
 	k_{\epsilon}^{(\alpha)}(x,y):=\frac{k_{\epsilon}(x,y)}{q_{\epsilon}^{\alpha}(x)q_{\epsilon}^{\alpha}(y)}\quad(\alpha\in\R).
\]
After carrying out the above procedure with the kernel $k=k_{\epsilon}^{(\alpha)}$ and the measure $\mu=q\,d\nu$, we obtain a Markov chain $P_{\epsilon,\alpha}$. Let $(\phi_l)_{l\in\N_0}\subset L_2(M,\nu)$ be the sequence of orthonormal eigenfunctions of the Laplace-Beltrami operator $\Delta$ on $M$ with Neumann boundary conditions. For $K\in\N_0$ let further
\[
	E_K:=\operatorname{span}\{\phi_0,\ldots, \phi_K \}.
 \]
It was derived in \cite[Theorem 2]{Coifman_Lafon_2006} that the spaces $E_K$ are contained in the domain of the generator of the Markov chain $P_{\epsilon,\alpha}$. More precisely, for such $f\in E_K$ it holds that
 \[
 	\lim_{\epsilon\to 0}\frac{P_{\epsilon,\alpha}-I}{\epsilon}f=\frac{\Delta(fq^{1-\alpha})}{q^{1-\alpha}}-\frac{\Delta(q^{1-\alpha})}{q^{1-\alpha}}f,
 \]
where the limit is taken in $L_2(M)$. Note that the authors of \cite{Coifman_Lafon_2006} use a different convention concerning the sign of the Laplacian: While we take $\Delta=\sum_{j=1}^n\partial_j^2$, the authors of \cite{Coifman_Lafon_2006} work with $\Delta=-\sum_{j=1}^n\partial_j^2$. The function $u(t):=P_{t,\alpha} u_0$ solves the equation
 \[
 	\partial_t u(t)=\frac{\Delta(u(t)q^{1-\alpha})}{q^{1-\alpha}}-\frac{\Delta(q^{1-\alpha})}{q^{1-\alpha}}u(t),\quad u(0)=u_0.
 \]
 Substituting $v(t):=q^{1-\alpha}u(t)$ leads to 
  \begin{align}\label{Eq:Diffusion_Maps:Schroedinger_TimeDependent}
 	\partial_t v(t)=\Delta v(t)-\frac{\Delta(q^{1-\alpha})}{q^{1-\alpha}}v(t),\quad u(0)=q^{1-\alpha}u_0.
 \end{align}
 Hence, one can study the operator
 \begin{align}\label{Eq:Diffusion_Maps:Schroedinger}
 	\psi\mapsto \Delta\psi -\frac{\Delta(q^{1-\alpha})}{q^{1-\alpha}}\psi
 \end{align}
and the dynamics generated by it in terms of the Markov chain $P_{\epsilon,\alpha}$. There are three values for $\alpha$ which are particularly interesting:
 \begin{itemize}
 	\item $\alpha=0$: This case is traditionally referred to in the literature as normalized graph Laplacian.
 	\item $\alpha=1$: In this case \eqref{Eq:Diffusion_Maps:Schroedinger} turns into the Laplace-Beltrami operator $\Delta$. In particular, the influence of the density $q$ is removed and the limit only depends on the geometry of the manifold.
 	\item $\alpha=\frac{1}{2}$: Here we obtain $\Delta-\frac{\Delta(q^{1-\alpha})}{q^{1-\alpha}}$. If the density is of the form $q=e^{-U}$, then \eqref{Eq:Diffusion_Maps:Schroedinger_TimeDependent} turns into 
 	\[
 		\partial_tv(t)=\Delta v(t)-\left(\frac{\|\nabla U\|_2^2}{4}-\frac{\Delta U}{2}\right)v(t)
 	\]
 	for \eqref{Eq:Diffusion_Maps:Schroedinger}, where $\|\nabla U\|_2$ denotes the Euclidean norm of the gradient of $U$. If we also substitute $f(t,x):=e^{-U(x)/2}v(t,x)$, then we obtain the equation
 	\[
 		\partial_t f(t,x)= \nabla\cdot(\nabla f(t,x)+f(t,x)\nabla U(x)),
 	\]
 	which is the forward Fokker-Planck equation to the stochastic differential equation
 	\begin{align}\label{Eq:Diffusion_Maps:SDE}
 		d X_t=-\nabla U(X_t)\,dt+\sqrt{2}\,dW_t
 	\end{align}
 	with reflecting boundary conditions at $\partial M$, where $(W_t)_{t\geq0}$ is a Brownian motion on $M$. Hence, the Markov chain $P_{\epsilon,\frac{1}{2}}$ can be used to study the probability densities of solutions of a stochastic differential equation of the form \eqref{Eq:Diffusion_Maps:SDE}.\\
 	In the discrete situation, the set $X$ is not given by a whole manifold $M$, but by a finite number of points sampled according to $q$. In this case integrals are replaced by sums, i.e. one works with
 	\begin{align*}
 		\bar{q}_{\epsilon}(x_i)=\sum_{j=1}^m k_{\epsilon}(x_i,x_j),\quad\bar{d}_{\epsilon}^{(\alpha)}(x_i)=\sum_{j=1}^m\frac{k_{\epsilon}(x_i,x_j)}{q_{\epsilon}(x_i)^{\alpha}q_{\epsilon}(x_j)^{\alpha},}\\
 		\bar{p}_{\epsilon,\alpha}(x_i,x_j)=\frac{k_{\epsilon}(x_i,x_j)}{\bar{d}_{\epsilon}^{(\alpha)}(x_i)},\quad \bar{P}_{\epsilon,\alpha}f(x_i)=\sum_{j=1}^m\bar{p}_{\epsilon,\alpha}(x_i,x_j)f(x_j)
 	\end{align*}
 where $m$ denotes the number of points. It was derived in \cite{Singer_2006} that 
 	\[
 		|P_{\epsilon,\alpha}f(x_i)-\bar{P}_{\epsilon,\alpha}f(x_i)|=O(m^{-\frac{1}{2}}\epsilon^{\frac{1}{2}-\frac{d_M}{4}})
 	\]
 	with high probability, where $d_M$ denotes the dimension of the underlying manifold. Therefore, the generators satisfy with high probability that
 	 	\[
 		|L_{\epsilon,\alpha}f(x_i)-\bar{L}_{\epsilon,\alpha}f(x_i)|=O(m^{-\frac{1}{2}}\epsilon^{-\frac{1}{2}-\frac{d_M}{4}}).
 	\]
 \end{itemize}
    
\smallskip
Diffusion maps have already found quite widespread applications across numerous areas~\cite{NadlerLafonCoifmanKevrekidis}
Of course, diffusion maps are just one instance of a broader class of data-driven reduction methods that are likely to gain more traction within climate science applications during the next decade.  In particular, we have not attempted to cover the full range of data-driven reduction methods such as in \cite{Kutz_2013} that presumably will be of great significance in future, including for climate model reduction; see for example \cite{Vlachas_2018,Scher_2018,Sapsis_2019}. Nonetheless, many of the explainable data-driven methods have their roots within reduction methods such as those discussed here, and indeed the methods discussed may contribute to such methods in future.

\section{Discussion}

In this review we have focused primarily on reduction methods of use (or potential use) in climate science that can be rigorously justified. This is a very diverse area and we unavoidably have had to limit ourselves to an incomplete set of reduction methods. 
These principles are likely to remain themes for future research as new techniques emerge and improve on existing reduction methods. Our brief review aims to make it easier to compare these methods: we try to highlight essential ideas for each method and provide references to concrete examples that serve as test cases. 

Probably the most relevant omission are reduction methods for stochastic differential equations, either SDEs or SPDEs, and more general data-based empirical reduction methods. We have not reviewed this area here, primarily to keep this work sufficiently accessible across disciplines, but there are undoubtedly many important methods also available for the stochastic case. As a basic principle, one can often generalise a reduction method for a deterministic differential equation to the stochastic case, usually at the expense of considerable additional technical work. 
Examples are path-wise stochastic slow manifold reductions via covariance tubes~\cite{BerglundGentz,KuehnSDEcont1}, reduction to Markov chain switching models via large deviation theory~\cite{FreidlinWentzell}, asymptotic analysis of stochastic slow manifolds~\cite{KabanovPergamenshchikov}, center/slow manifold reduction using a skew-product flow formalism~\cite{Boxler,DuanLuSchmalfuss,KuehnNeamtu1,SchmalfussSchneider}, distribution-based approaches simplifying Fokker-Planck equations~\cite{Schuss}, nonlocal Mori-Zwanzig type reductions~\cite{ChekrounLiuWang,ChekrounLiuWang1}, moment closure~\cite{Kuehn_2016} and amplitude/modulation equations~\cite{Bloemker}, just to name a few. 

One might conjecture that a more coherent and rigorous view of the area will emerge within the coming years as some reduction methods for stochastic systems have been applied successfully in practical climate models already using numerical approaches~\cite{DijkstraEtal:2016}. It is clear that such reductions may be very informative, especially statistical mechanics and large deviations approaches to weather/climate extreme problems; see for example \cite{bouchet2012statistical,galfi2019large,galfi2021fingerprinting,margazoglou2021dynamical} as well as for multistable regimes of the climate system \cite{lucarini2017edge}. Another area where rigorous reduction methods are challenged is in justifying the predominantly empirical sub-grid parametrizations used in climate models \cite{majda1999models,chorin1998optimal,palmer2005representing,kwasniok2012data}.

It is important to acknowledge that many reduction methods of great utility for climate models are being, and need to be, applied in cases where it is impossible to provide rigorous proofs of the reduction along the lines presented here. Nonetheless, rigorous reduction methods are still useful in such contexts because better understanding of when and whether these methods are appropriate. In particular, these can (a) help to find an optimal set of conditions and results such that reduction method can be applied and (b) help to understand the circumstances and details of when such a reduction method fails to work, and what happens in cases where it cannot be applied.

\subsection*{Acknowledgements}

This project received funding from the European Union's Horizon 2020 research and innovation programme under grant agreement No 820970 (TiPES). CK also would like to acknowledge support by the VolkswagenStiftung for support via a Lichtenberg Professorship.

\bibliographystyle{abbrv}
\bibliography{Bibliography.bib}

\end{document}